\documentclass[final]{siamltex}
\usepackage{amsmath}
\usepackage{amssymb}
\usepackage{epsfig}

\newtheorem{defi}{Definition}[section]
\newcommand{\bd}{\begin{defi}} 
\newcommand{\ed}{\end{defi}}

\newcommand{\bl}{\begin{lemma}}
\newcommand{\el}{\end{lemma}} 

\newcommand{\bt}{\begin{theorem}}
\newcommand{\et}{\end{theorem}}
	
\newtheorem{remark}[defi]{Remark}
\newcommand{\br}{\begin{remark}}
\newcommand{\er}{\end{remark}}

\newtheorem{cor}[defi]{Corollary}
\newcommand{\bc}{\begin{cor}}
\newcommand{\ec}{\end{cor}}

\newtheorem{pro}[defi]{Proposition}
\newcommand{\bp}{\begin{pro}}
\newcommand{\ep}{\end{pro}}

\def\R{\mathbf{R}}
\def\C{\mathbf{C}}
\def\bs#1{\boldsymbol{#1}}
\def\lra#1{\langle{#1}\rangle^{}}

\def\Uhb{U_{h,b}^{}}
\def\Uh{U_h^{}}
\def\pUh{\tilde{U}_h^{}}
\def\pUhm{\tilde{U}_{h,-}^{}}
\def\Uhi{U_{h,i}^{}}

\def\Uhm{U_{h,-}^{}}
\def\Uhmn{U_{h,-}^{n}}
\def\Uhmnn{U_{h,-}^{n+1}}

\def\Uhp{U_{h,+}^{}} 
\def\Uhpn{U_{h,+}^{n}} 
\def\Uhpnn{U_{h,+}^{n+1}}

\def\Uhpm{U_{h,\pm}^{}}

\def\Uhpmnn{U_{h,\pm}^{n+1}}
\def\Uhmpn{U_{h,\mp}^{n}}

\def\U0{U_{0}^{}}
\def\nb{\nabla}
\def\nh{\nabla_h^{}}
\def\dh{\nabla_h^{}\cdot}
\def\Ubh{{\overline{U}}_h^{}}
\def\Ubhm{{\overline{U}}_{h,-}^{}}

\def\Ubhp{{\overline{U}}_{h,+}^{}}

\def\La{\Delta}
\def\Lah{\Delta_h^{}}
\def\Om{\vec{\Omega}}
\def\ub{\overline{u}}

\def\ubmnn{\overline{u}^{n+1}_-}

\def\vb{\overline{v}}

\def\ha{\hat{\zeta}}

\def\hU{\hat{U}_h}   

\def\hUo{\hat{U}_h^0} 
\def\huo{\hat{u}^0}
\def\hvo{\hat{v}^0}
\def\hUn{\hat{U}_h^n}   
\def\hUnm{\hat{U}_{h,-}^{n}}   
\def\hUnp{\hat{U}_{h,+}^{n}}   
\def\hUnmp{\hat{U}_{h,\mp}^{n}}   

\def\hXm{\hat{X}_{-}} 
\def\hXp{\hat{X}_{+}}

\def\en{\bs{e}_n}
\def\a{\zeta}
\def\pa{\tilde{\zeta}}
\def\pam{\tilde{\zeta_-}}
\def\ai{\zeta_i}
\def\ab{\zeta_b}
\def\am{\zeta_-}

\def\amnn{\zeta_-^{n+1}}

\def\ap{\zeta_+}

\def\apnn{\zeta_+^{n+1}}

\def\apm{\zeta_\pm}

\def\apmnn{\zeta_\pm^{n+1}}

\def\b{\chi}
\def\eps{\varepsilon}
\def\seps{\sqrt{\varepsilon}}
\def\dx{\Delta x}

\def\dz{\Delta z}
\def\dt{\Delta t}

\def\pt{\partial_t}
\def\px{\partial_x}
\def\py{\partial_y}
\def\pz{\partial_z}
\def\pn{\partial_n}
\def\Dx{D_x}

\def\Dz{D_z}

\def\sxx{\sigma_{xx}}
\def\sxy{\sigma_{xy}}
\def\sxz{\sigma_{xz}}
\def\szz{\sigma_{zz}}
\def\szx{\sigma_{zx}}
\def\szy{\sigma_{zy}}
\def\syx{\sigma_{yx}}
\def\syy{\sigma_{yy}}
\def\syz{\sigma_{yz}}

\def\eps{\varepsilon}

\DeclareMathOperator{\EP}{PE}

\title{Optimized Schwarz waveform relaxation for Primitive Equations of the ocean}

\author{ E. Audusse, P. Dreyfuss, B. Merlet. 
\thanks{Universit\'e Paris Nord - Institut Galil\'ee
LAGA (Laboratoire d'Analyse, G\'eom\'etrie et Applications)
Avenue J.B. Cl\'ement, 93430 Villetaneuse (merlet@math.univ-paris13.fr).}}

\begin{document}
\maketitle

\begin{abstract}
In this article we are interested in the derivation 
of efficient domain decomposition methods 
for the viscous primitive equations of the ocean. 
We consider the rotating 3d incompressible 
hydrostatic Navier-Stokes equations with free surface. 
Performing an asymptotic analysis of the system 
with respect to the Rossby number, we compute an approximated
Dirichlet to Neumann operator and build an optimized Schwarz waveform relaxation algorithm.
We establish the well-posedness of this algorithm and present some numerical results to illustrate the method.
\end{abstract}



\pagestyle{myheadings}
\thispagestyle{plain}
\markboth{}{}

\section{Introduction}
\label{secintro}
A precise knowledge of ocean parameters (velocity, temperature...)
is an essential tool to obtain climate and meteorological previsions.
This task is nowadays of major importance 
and the need of global or regional simulations
of the evolution of the ocean is strong. 
Moreover the large size of global simulations 
and the interaction between global and regional models 
require the introduction of efficient domain decomposition methods.\\
\newline
The evolution of the ocean is commonly modelized by the use of the viscous primitive equations. 
This system is deduced from the full three dimensional incompressible Navier-Stokes equations 
with free surface with the use of the hydrostatic approximation and of the Boussinesq hypothesis. 
It is implemented in all the major softwares that are concerned 
with global or/and regional simulations of ocean and/or atmosphere
(we refer for example to NEMO \cite{nemo}, MOM \cite{mom} or HYCOM for global models 
and ROMS \cite{roms} or MARS for regional models).
The primitive equations have been studied for twenty years 
and important theoretical results are now available \cite{lionspere,temam04,titi}. 
The numerical treatment of this system has been also strongly investigated \cite{temam08}.
But the key point here is to simulate global circulation on the earth 
for long time and/or with small space discretization. 
This type of computations can not be performed 
on a single computer in realistic CPU time and need to be parallelized. 
The problem is then to allow the different
subdomains to interact in an efficient way. 
Another type of applications that are commonly investigated
in the oceaonographic and/or meteorological community
is to couple global and regional models in order to obtain
precise regional previsions. The problem 
is also to construct an efficient interaction between the two models.
In \cite{blayo} the authors exhibit that most of the existing algorithms 
are not able to compute this kind of problem in an efficient way.
We propose in this article to investigate these still open questions 
in the context of a quite recent performing domain decomposition method : 
the Schwarz waveform relaxation type algorithms.\\
\newline
The development of domain decomposition techniques 
have known a great development for the last decades
and our purpose is not to make an exhaustive 
presentation of these methods. We refer the reader to \cite{quarteroni,toselli} for a general presentation
and we restrict ourselves to the description of Schwarz waveform relaxation method. 
It is a relatively new domain decomposition technique.
It has been developed for the last decade 
and has been successfully applied to different types of equations.
This type of algorithms is the result of the interaction
between classical Schwarz domain decomposition techniques and waveform relaxation algorithms.
Its great interest is to be explicitly designed for evolution equations
and to allow different strategies for the space time discretization 
in each subdomain. Moreover we can even consider different models
in each subdomain without modifying the architecture of the interaction.\\
The heart of the classical Schwarz method is to solve 
the problem on the whole domain thanks to an iterative procedure
where a problem is solved on each subdomain by the use 
of boundary conditions that contain the information coming from the neighboring subdomains.
It comes from the early work of Schwarz \cite{schwarz} where this idea was introduced
to prove the well-posedness of a Poisson problem in some nontrivial domains.
This method is designed for stationary problems and presents 
two main drawbacks : it needs an overlapping between subdomains and it
converges slowly \cite{lions1}.
In the last decade, some works have been devoted to cure these disagreements \cite{lions3}.
We refer to \cite{gander06} for a complete presentation.\\
The extension to time evolution problems was performed at the end of the nineties
by Gander \cite{gander97,gander98} and Giladi \& Keller \cite{giladi} 
and was denoted Schwarz waveform relaxation algorithms.
The authors mixed the classical Schwarz approach 
with waveform relaxation techniques
developed in the context of the solutions of large system 
of ordinary differential equations \cite{lelarasmee,jeltsch}.
The exchanged quantities were of Dirichlet type.
Optimized Schwarz waveform relaxation methods were developed
with the introduction of more sophisticated information 
to compute the interaction between the subdomains.
These optimized algorithms were based on previous works \cite{engquist,halpern86,halpern91} 
about the derivation of absorbing boundary conditions 
respectively for hyperbolic, elliptic and incompletely parabolic equations.
The same ideas were used to derive efficient transmission conditions
between the subdomains : since the exact transparent conditions 
can not be implemented in general
(it may lead to non-local pseudo-differential operators), 
the derivation of some approximate conditions is performed. 
These conditions can be optimized with respect to some free parameters 
which justifies the name of the method.
The optimized Schwarz waveform relaxation method
was first applied to the wave equation \cite{halpern03} and then to the advection-diffusion equation with 
constant or variable coefficients~\cite{martin05}.
A recent paper \cite{halpern07} gives the complete solution of the one dimensional optimization problem 
for constant coefficients equations.
More recently the method has been extended 
to the linearized viscous shallow water equations 
without advection term by V. Martin~\cite{martin08}. 
Here we are interested in the application 
of the method to the system of Primitive Equations of the ocean. 
It leads to non-trivial new problems 
(new transmission conditions, well-posedness of the problem, convergence of the algorithm...)
that we address in this article.\\ 
\newline
The outline of the paper is the following : in Section~\ref{secequations} we write the equations
and we precise the asymptotic regime that we consider.
In Section~\ref{secalgo} we derive an approximated Dirichlet to
Neumann operator, and define the associated Schwarz waveform
relaxation algorithm. In 
Section~\ref{secwp} we define a weak formulation of the problem on
the whole domain and prove that it is well-posed in the natural
functional spaces. In Section~\ref{secWP} we introduce a weak
formulation for the Schwarz waveform relaxation algorithm and prove that each sub-problem 
solved in the algorithm is well-posed.
Finally we present some numerical results 
in Sections~\ref{secnum}. 
\section{The set of equations}
\label{secequations}
We first write the primitive equations of the ocean.
Then we present the simplified system from which we are able 
to derive efficient transmission conditions.
\subsection{The primitive equations of the ocean}
We consider the primitive equations of the ocean 
on the domain $(x,y,z,t) \in \R \times \R \times [-H(x,y),\a(x,y,t)] \times \R^+$
where $-H(x,y)$ denotes the topography of the ocean 
and $\a(x,y,t)$ denotes the altitude of the free surface of the ocean.
The primitive equations are commonly written~\cite{C-R94}
\begin{eqnarray}
\label{eqprim}
\pt \Uh +\Uh\cdot \nh \Uh -\nu \La \Uh + \frac{2}{\rho_0}\Om\wedge \Uh +\cfrac1{\rho_0} \nh p &=&0,\\
\label{div}
\dh \Uh + \pz w &=&0,\\
\label{p}
\pz p &=& -\rho g,\\
\label{rho}
\rho &=& \rho(z,T,S),\\
\label{T}
\pt T +\U0\cdot \nb T -\nu_T \La T  &=& Q_T,\\
\label{S}
\pt S +\U0\cdot \nb S -\nu_S \La S  &=& Q_S,
\end{eqnarray}
where the unknowns are the 3d-velocity $(\Uh,w)=(u,v,w)$, the pressure $p$,
the density $\rho$, the temperature $T$ and the salinity $S$.
The parameters are the gravity $g$, the eddy viscosity $\nu$,
the eddy diffusion coefficients for the tracers $\nu_T$ and $\nu_S$
and the earth rotation vector $\Om$. The source terms $Q_T$ and $Q_S$ for the 
temperature and salinity model the influence of the sun, rivers and atmosphere 
for these tracers.\\ 
Note that we consider here the classical but non-symmetric viscosity tensor
\begin{eqnarray*}
\sigma &=&
\left(\begin{array}{ccc}
\sxx& \sxy & \sxz\\
\syx& \syy & \syz\\
\szx & \szy & \szz 
\end{array}\right)
\ =\  \nu
\left(\begin{array}{ccc}
 \px u & \py u & \pz u \\
 \px v & \py v & \pz v \\ 
 \px w & \py w & \pz w \\
\end{array}\right).
\end{eqnarray*}
Other form of the viscosity tensor can be found in \cite{perthame}.
Note also that it is possible to consider different viscosity coefficients 
in the horizontal and vertical directions \cite{rousseau}.\\
These equations are supplemented by initial and boundary conditions.
At initial time, we impose 
\begin{equation*}
\Uh(\cdot,0)\ =\ \Uhi\quad\mbox{in } \Omega,\quad  \qquad
\a(\cdot,0)=\ai\quad\mbox{in }\omega,
\end{equation*}
where the subscript letters $i$ means ``initial''. At the bottom of the ocean we impose a non-penetration condition and 
a friction law of Robin type ($\alpha_b > 0$)
\begin{equation}
\label{bottom}
\Uh(-H) \cdot \nh (H) - w(-H)=0, \qquad \pn U_t (-H) + \alpha_b U_t (-H) = 0,
\end{equation}
where $U_t$ stands for the tangential velocity
and $n$ denotes the outward normal vector to the bottom of the ocean.\\
The free surface is transported by a kinematic boundary condition
\begin {equation}
\label{CCSL}
\pt \a + \Uh(\a) \cdot\nh \a -w(\a)\ =\ 0.
\end{equation}
The equilibrium of the stresses at the free surface implies 
\begin{equation}
\label{tensions}
\left[\sigma - (p-p_a) Id\right] 
\cdot
\frac{1}{\sqrt{1+(\px \a)^2+(\py \a)^2}}
\left(\begin{array}{c}
\px \a\\
\py \a\\
1
\end{array}
\right)\ =\ 0,
\end{equation}
where $p_a(x,y,t)$ denotes the atmospheric pressure.\\
\begin{figure}[h]
\centering
\input{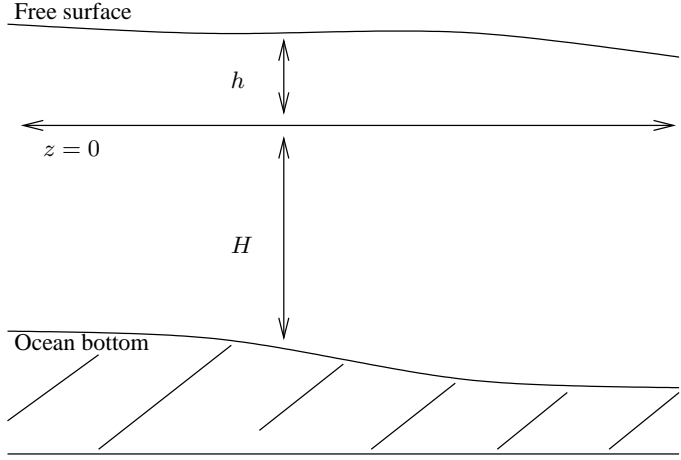} 
\caption{\label{ocean}Schematic representation of the ocean}
\end{figure}

\subsection{A linearized hydrostatic model}
In order to derive simple and efficient transmission conditions
for the Schwarz waveform relaxation method we make 
some assumptions on this set of equations. \\
\newline
First we neglect the influence of the tracers (temperature and salinity)
on the density. Thus we suppose that the density is constant (we assume $\rho_0=1$) 
and we do not solve the equations on the tracers~\eqref{T}-\eqref{S}. 
Note that these equations are classical advection-diffusion equations 
for which the optimized transmission conditions are well known
\cite{halpern07,martin05}.\\
Then we use the divergence-free condition~\eqref{div} and the non-penetration condition~\eqref{bottom}
to write the vertical velocity $w$ as a function of the horizontal velocity $\Uh$
and we use the hydrostatic assumption~\eqref{p}
to write the pressure $p$ as a function of the water height $\a$.
The remaining unknowns in the system are the horizontal velocity $\Uh$ and the water height $\a$.
The set of equations~\eqref{eqprim}-\eqref{p} stands
\begin{eqnarray*}
\pt \Uh +\Uh\cdot \nh \Uh -\nu \La \Uh + f C \Uh + g \nh \a &=&0,\\
\pt \a + \nh \cdot \int_{-H}^\a \Uh \, dz &=&0,
\end{eqnarray*}
with 
\begin{eqnarray*}
C=
\left(
\begin{array}{cc}
0 & -1\\
1  & 0
\end{array}
\right),&\qquad\mbox{and}\quad  f:=2\Om\cdot\bs{e}_z.
\end{eqnarray*}
The first equation is written on the initial domain $\R_x \times \R_y \times [-H(x,y),\a(x,y,t)]_z \times \R^+_t$
while the second one is written on $\R_x \times \R_y \times \R^+_t$.
We consider for simplicity a flat bottom
and a constant atmospheric pressure.
Then we linearize the problem around a constant state 
which corresponds to a horizontal velocity $\U0=(u_0,v_0)$ and a horizontal free surface located at $z=0$.
It follows that the water height $\a$ is a small perturbation.
In the sequel $\Uh$ denotes the perturbation on the horizontal velocity.
The linearized problem stands
\begin{eqnarray}
\label{eqprimlin}
\pt \Uh +\U0 \cdot \nh \Uh -\nu \La \Uh + fC \Uh + g \nh \a &=&0,\\
\label{divlin}
\pt \a + H \nh \cdot \Ubh + \U0 \cdot \nh \a &=&0,
\end{eqnarray}
where $$\Ubh \ =\ \left(
\begin{array}{c}
\overline{u} \\
\overline{v} 
\end{array}
\right) := \frac{1}{H} \int_{-H}^0 \Uh \,dz,$$ denotes the mean horizontal
velocity of the flow. This mean velocity is called barotropic velocity by the oceanographic 
community while the deviation $\Uh-\Ubh$ is called baroclinic velocity~\cite{C-R94}.\\
Note that the first equation~\eqref{eqprimlin} is now written in the fixed domain 
$\R_x \times \R_y \times [-H,0]_z \times \R^+_t$.\\
The associated boundary conditions are
\begin{eqnarray}
\label{boundlin}
\pz \Uh (z=0) \ =\  0, \qquad \left\{\pz \Uh  + \alpha_b \Uh\right\}
(z=-H) \ =\ 0,
\end{eqnarray}
where the boundary condition at $z=0$ is deduced from the equilibrium of the stresses at the free surface~\eqref{tensions}.
Indeed, with the help of the linearization procedure, we first deduce
that at first order $\pz \Uh(z=\a) \simeq 0$ and then $\pz \Uh(z=0) \simeq 0$ since we
assume that $\a$ is small.\\
In order to derive the transmission conditions we assume $\alpha_b = 0$ in the sequel. 
However in the definition of the Schwarz waveform relaxation algorithm and in the numerical simulations,
the condition $\alpha_b > 0$ will be supported. 
\subsection{Dimensionless system}
We choose characteristic horizontal and vertical lengths 
(denoted $L$ and $H$ respectively) 
and velocity $U$ of the problem.
We introduce the dimensionless quantities
\begin{eqnarray*}
(x,y)= L(\tilde{x},\tilde{y}),    &\qquad & t = (L/U) \tilde{t},\\
\a =H \tilde{\a}, &\qquad & z =H \tilde{z},\\
\Uh=U \widetilde{\Uh}, & \qquad & \U0 = U \widetilde{\U0}.   
\end{eqnarray*}
The spatial domains of computation are $\Omega = \R_x \times \R_y \times (-1,0)_z$
for the momentum equation and $\omega = \R_x \times \R_y$ for the
continuity equation. We study both equations in the time interval
$[0,T]$, where $T > 0$ is fixed. \\
Dropping the ``$\sim$'' for a better readability, the system in dimensionless variables stands 
\begin{eqnarray}
\label{eqprimlinadim}
\partial_t \Uh +\U0\cdot \nh \Uh -\cfrac1{Re} \Lah \Uh -\cfrac1{Re'}
\pz^{\,2} \Uh +\cfrac1\eps  C\Uh +\cfrac1{Fr^2} \nh\a&=&0,
\\
\label{boundlinadim}
\pz \Uh (x,y,0,t)=\pz \Uh (x,y,-1,t)&=&0,\\
\label{initu}
\Uh(\cdot,0)&=&\Uhi,\\
\label{divlinadim}
\partial_t \a +\U0 \cdot \nh \a +\dh \Ubh &=&0,\\
\label{inita}
\a(\cdot,0)&=&\ai.
\end{eqnarray}
We have introduced the characteristic quantities,
\begin{description}
\item[-] $\eps=U/(fL)$ the Rossby number, 
\item[-]$Re=UL/\nu$ the horizontal Reynolds number, 
\item[-] $Re'=H^2/L^2 Re$ the vertical Reynolds number,
\item[-]  $Fr=U/\sqrt{gH}$ the Froude number. \\
\end{description}
We choose to exhibit the Rossby number as a small parameter
since we are interested in long-time oceanographic circulation
for which the Rossby number is typically of magnitude $10^{-2}$.
The values of Reynolds and Froude numbers vary with respect 
to the turbulent processes and to the depth of the area that is considered respectively.
\section{The optimized Schwarz waveform relaxation algorithm}
\label{secalgo}
We are now interested in finding efficient transmission conditions 
for equations~\eqref{eqprimlinadim}--\eqref{inita}.
We first present the Schwarz waveform relaxation method.
Then we derive the relations satisfied by the {\it optimal}
transmission conditions. Since we are not able to solve analytically these
equations, we perform an asymptotic analysis with respect to the Rossby number $\eps$
in order to derive some {\it approximated} transmission conditions.
Finally we present the related {\it optimized} Schwarz waveform relaxation algorithm.
\subsection{The Schwarz waveform relaxation method}
The heart of the method is the following.  We first divide the computational domain into an arbitrary number of subdomains.
Then we solve each sub-problem independently for the whole time interval. 
The interactions between neighboring subdomains are entirely contained in the boundary conditions.
An iterative procedure is considered until a prescribed precision is reached.
The advantages of the method are clear : the parallelization is almost optimal: at each step 
the sub-problems are solved independently, so the space-time discretization strategies
(or even the models...) can be chosen independently on each subdomain. Moreover at the end of each step only a small amount of informations are exchanged.  
The main drawback is related to the needed number of iterations : 
the method is efficient if it converges quickly (in two or three iterations typically).
This requirement needs the derivation of efficient transmission
conditions. \\
In the sequel we consider for simplicity two subdomains
but the method extends to an arbitrary number of
subdomains.\\ 
\newline
We begin with some notations. First we introduce the left and right spatial subdomains $\Omega^-$ and $\Omega^+$ defined by:  
\begin{eqnarray*}
\Omega^- \ :=\  (-\infty,0)_x \times \R_y \times (-1,0)_z,& 
\qquad &\Omega^+ \ :=\  (0,+\infty)_x \times \R_y \times (-1,0)_z, 
\end{eqnarray*}
and their interface  
\begin{eqnarray*}
\Gamma &= &\{0\}_x \times \R_y \times (-1,0)_z\ \simeq\ \R_y\times(-1,0)_z.
\end{eqnarray*}
We also introduce the domains $\omega^\pm := \pm (0,+\infty)_x \times \R_y$
for the unknowns that do not depend on the $z$ variable, 
and their interface $\gamma := \{0\}_x \times \R_y\simeq \R_y$.\\
Let $D$ be some spatial open domain and $T>0$ be a given real
number. Then we will write $D_T$ to denote the cylindrical domain $D_T := D
\times (0,T)$.\\
We denote by $\EP$ the set of equations 
\eqref{eqprimlinadim},~\eqref{boundlinadim} and~\eqref{divlinadim} and
$X:=(\Uh,\a)$ stands for the solution of this system with associated
initial data $X_i:=(\Uhi,\ai)$.\\
Then the Schwarz waveform relaxation algorithm is defined as follows: 
\begin{equation}
\left\{
\begin{array}{rll} \EP(X^{n+1}_-) =&  0&
  \mbox{on }  \Omega^-_T, \\
&&\\
X^{n+1}_-(\cdot,0) =& X_i&\mbox{on } \Omega^- , \\
&&\\
\mathcal{B}_- X^{n+1}_- =& \mathcal{B}_- X^n_+ &\mbox{on } \Gamma_T, 
\end{array}
          \right. 
                   ~\qquad
\left\{
\begin{array}{rll} 
\EP(X^{n+1}_+) =&  0  &\mbox{on } \Omega^+_T\\
 &&\\
X^{n+1}_+(\cdot,0) =& X_i&\mbox{on } \Omega^+,  \\
&&\\
\mathcal{B}_+ X^{n+1}_+ =& \mathcal{B}_+ X^n_- &\mbox{on }  \Gamma_T, 
\end{array}
           \right.~
\label{swra}
\end{equation}
where the operators $\mathcal{B}_\pm$ 
contain the transmission conditions.\\
In the classical Schwarz waveform relaxation algorithm \cite{gander97,gander01},
the transmitted quantities are of Dirichlet type 
and the operators $\mathcal{B}_\pm$ are thus chosen to be the identity operator.
Note that in this case an overlap is needed in the definition of the subdomains.\\
In the sequel we are interested in deriving more efficient transmission conditions.
In order to reach such a goal 
we will first describe the general method 
to obtain optimal transmission conditions.
The transmission conditions are said to be optimal 
if the algorithm converges in two iterations to the solution of the initial problem.
These optimal transmission conditions involve the Dirichlet to Neumann operator associated to $\EP$ on the subdomains $\Omega^\pm_T$.  Here we will see that we are not able to obtain an explicit
formulation for these optimal conditions. Anyway these optimal
boundary conditions are not local and consequently too expensive to be
useful from a numerical point of view.\\
Recent methods have been developed recently in order to approximate 
these optimal conditions by analytical or numerical
means --- see the review paper~\cite{gander06} for elliptic problems and~\cite{halpern07} for parabolic evolution equations.\\
Here we will perform an asymptotic analysis of the system  with
respect to the Rossby number $\eps$ in order to deduce 
a set of approximated and efficient transmission conditions. This strategy has been initiated in~\cite{martin08} for the shallow water equation without advection term. 
In our case, it turns out that these approximate transmission conditions lie in a 
two parameter family of boundary conditions. In Section~\ref{secnum} we
optimize numerically  the transmission conditions in this two parameter
family.\\
\newline
Let us first describe in a formal setting the ideal case of optimal transmission conditions for the Schwarz waveform relaxation algorithm~\eqref{swra}.\\
We consider the case $u_0 > 0$. 
The case $u_0 < 0$ is deduced by applying the symmetry ``$x'=-x$''.
Integrating the linearized Primitive equations~$\EP$ on a subdomain, we see that the flux of the unknown $(\Uh,\a)$ through the interface $\Gamma$ is given by
$$
\left(\displaystyle{\frac{1}{Re}} \px \Uh - u_0 \Uh 
- 
\displaystyle{\frac{1}{Fr^2}} 
\left(
\begin{array}{c}
\a\\
0
\end{array}
\right) \,,\, 
u_0 \a + {\overline{u}}
\right)
$$
Using this flux as a Neumann operator, we define the Dirichlet to Neumann operators as follows. Consider a Dirichlet data $X_b=(\Uhb,\a_b)$, we set 
$$DN_-^{\Uh} X_b\ :=\ \left(\displaystyle{\frac{1}{Re}} \px \Uh - u_0 \Uh 
- 
\displaystyle{\frac{1}{Fr^2}} 
\left(
\begin{array}{c}
\a\\
0
\end{array}{c}
\right)\right)_{|\Gamma_T},$$ 
where $X=(\Uh,\a)$ solves
\begin{equation*}
\left\{
\begin{array}{rll} 
\EP(X) =&  0  &\mbox{on } \Omega^+_T\\
X(\cdot,0) =& 0&\mbox{on } \Omega^+,  \\
X =& X_b &\mbox{on }  \Gamma_T.
\end{array}
           \right.~
\end{equation*}
Symmetrically, consider a Dirichlet data $\Uhb$, we set 
$$\left(\begin{array}{c}DN_+^{\Uh} X_b\\DN_+^{\a}  X_b\end{array}\right) \ :=\ \left(\begin{array}{l} \displaystyle{\frac{1}{Re}} \px \Uh - u_0 \Uh 
- 
\displaystyle{\frac{1}{Fr^2}} 
\left(
\begin{array}{c}
\a\\
0
\end{array}
\right)\\
u_0 \a + {\overline{u}}
\end{array}\right)_{|\Gamma_T},$$ 
where $X=(\Uh,\a)$ solves
\begin{equation*}
\left\{
\begin{array}{rll} \EP(X) =&  0&
  \mbox{on }  \Omega^-_T, \\
X(\cdot,0) =& 0&\mbox{on } \Omega^- , \\
\Uh =& \Ubh &\mbox{on } \Gamma_T.
\end{array}
\right. 
\end{equation*}
Notice that since we consider the case $u_0 > 0$, the continuity equation~\eqref{divlinadim}
the boundary condition is relevant only in the subdomain $\Omega^+_T$. This is why in the later case we do not have to prescribe a boundary condition for $\a$.\\

\noindent
Once these Dirichlet to Neumann operators are defined we can introduce the optimal transmission conditions 
\begin{eqnarray}
\mathcal{B}_- X &=& 
\left(
\begin{array}{c}
\displaystyle{\frac{1}{Re}} \px \Uh - u_0 \Uh 
- 
\displaystyle{\frac{1}{Fr^2}} 
\left(
\begin{array}{c}
\a\\
0
\end{array}
\right)
- DN^{\Uh}_- X \\
\end{array}
\right),
\label{algoopt-}
\\
\mathcal{B}_+ X &=& 
\left(
\begin{array}{c}
-\displaystyle{\frac{1}{Re}} \px \Uh + u_0 \Uh + 
\displaystyle{\frac{1}{Fr^2}} 
\left(
\begin{array}{c}
\a\\
0
\end{array}
\right)
- DN^{\Uh}_+ X\\
\\
u_0 \a + {\overline{u}}- DN^\a_+ X
\end{array}
\right),
\label{algoopt+}
\end{eqnarray}

\medskip

\bp
With this particular choice of transmission operators $\mathcal{B}_\pm$, the algorithm~\eqref{swra} converges in two iterations.\\  
\ep

\noindent
\begin{proof}
By linearity, we may assume that the exact solution is $0$ ($X_i\equiv 0$).
At the initial step the solutions on each subdomain  do not satisfy any particular property.
But the first iterate solves the primitive equations with vanishing initial data.
It follows from the very definition of the operators $DN_\pm$ that in the
definition of the second iterate, the right hand sides of the
transmission conditions vanish for both
sub-problems. We deduce that this second iterate vanish: the algorithm converges in two steps.
\end{proof}

\medskip

\noindent
The operators~\eqref{algoopt-}\eqref{algoopt+} being non-local pseudo-differential operator, they are not well suited for numerical implementation. Our strategy is to approximate these operators by numerically cheap operators. Of course the two-step convergence property will be lost. The quality of the approxiamation will be measured through the convergence rate of the algorithm. From the structure of~\eqref{algoopt-}\eqref{algoopt+}, we choose to write $\mathcal{B}_\pm$ as perturbations of the natural operators transmitted through the interface:
\begin{eqnarray}
\mathcal{B}_- X &=& 
\left(
\begin{array}{c}
\displaystyle{\frac{1}{Re}} \px \Uh - u_0 \Uh 
- 
\displaystyle{\frac{1}{Fr^2}} 
\left(
\begin{array}{c}
\a\\
0
\end{array}
\right)
- \mathcal{S}^{\Uh}_- X \\
\end{array}
\right),
\label{algooptimal-}
\\
\mathcal{B}_+ X &=& 
\left(
\begin{array}{c}
-\displaystyle{\frac{1}{Re}} \px \Uh + u_0 \Uh + 
\displaystyle{\frac{1}{Fr^2}} 
\left(
\begin{array}{c}
\a\\
0
\end{array}
\right)
- \mathcal{S}^{\Uh}_+ X\\
\\
u_0 \a + {\overline{u}}- \mathcal{S}^\a_+ X
\end{array}
\right),
\label{algooptimal+}
\end{eqnarray}
where $\mathcal{S}^{\Uh}_\pm$ and $\mathcal{S}^\a_+$ are
pseudo-differential operators that will approximate the Dirichlet to Neumann operators.\\
Let us finally remark that the differences in the expression of 
the two transmission operators $\mathcal{B}_\pm$ are due to the sign $u_0>0$.
Since $\mathcal{B}_-$ contains the information
that is transmitted from $\Omega^+_T$ to $\Omega^-_T$
it is constructed on three boundary values (velocities and water height)
but it has to transmit only two boundary conditions for momentum equations~\eqref{eqprimlinadim}.
On the contrary $\mathcal{B}_+$ is constructed on two boundary values (velocities) 
but has to send three boundary conditions (for momentum and continuity equations).
\\
In the next subsections we will identify optimal and approximated transmission operators. To carry out the computation of the Dirichlet to Neumann operators we perform Fourier-Laplace transforms.

\subsection{Laplace-Fourier transform of the primitive equations}
We perform on the set of primitive equations~\eqref{eqprimlinadim}-\eqref{divlinadim}
a Fourier transform in the $y$ variable and a Laplace transform in time.
The dual variables are respectively denoted $\eta \in \R$ and $s= \sigma + i \tau
\in \C$. The real part $\sigma$ is assumed to be strictly positive.
We obtain in each subdomain the same set of differential equations
\begin{eqnarray*}
\left\{ s  + u_0 \px  + i \eta v_0  - \cfrac1{Re} {\px^{\,2}}  
+ \cfrac1{Re} \eta^2  - \cfrac1{Re'} {\pz^{\,2}}  
+ \cfrac1\eps B\right\} \hU + \cfrac1{Fr^2} 
\left(
\begin{array}{c} 
\px\\
i\eta 
\end{array}
\right)
\ha
&=&0,
\nonumber\\
\\
\left\{s  + u_0 \px  + i \eta v_0 \right\} \ha + \px {\hat{\ub}} + i \eta {\hat{\vb}} &=&0.\nonumber\\
\end{eqnarray*}
In the $z$ direction we introduce the eigenmodes of the operator $- {\pz^{\,2}}$
on $(-1,0)$ with homogeneous Neumann boundary conditions~\eqref{boundlinadim}
\begin{eqnarray*}
\en(z) &:=& \alpha_n\cos(\mu_n z) 
\quad \mbox{with} \quad \mu_n :={n\pi}; 
\quad \alpha_0:=1 \quad \mbox{and} \quad \alpha_n:={\sqrt{2}} \quad \mbox{if} \quad n>0. 
\end{eqnarray*}
Then we search for the solution on the form
\begin{eqnarray*}
\hU(x,z)= \sum_{n=0}^\infty \hUn(x) e_n(z).
\end{eqnarray*}
Note that we obviously obtain
${\hat{\Ubh}} = \hUo.$
It means that the first vertical mode $\hUo$ represents the barotropic velocity
while the sum of the other ones denotes the baroclinic deviation.\\
The barotropic mode is coupled with the water height and 
it is the solution of the following system of three ordinary differential
equations, 
\begin{eqnarray}
\label{eqprimlinadim0}
-\cfrac1{Re} \px^{\,2}
\hUo  +u_0\partial_x \hUo + \left\{ s  + i\eta  v_0 + \cfrac1{Re}
\eta^2  +\cfrac1\eps  B\right\}\hUo +\cfrac1{Fr^2}
\left(\begin{array}{c} \partial_x\ha \\ 
i\eta \ha\end{array} \right) = 0,
\\
\label{divlinadim0}
u_0\partial_x \ha +\left(s+ i\eta v_0\right) \ha 
+\partial_x \huo+i\eta \hvo = 0.
\end{eqnarray}
This last system is exactly the Laplace-Fourier transform 
of the so-called linearized viscous shallow water equations
\cite{stvenant}. \\ 
\newline
For the other vertical modes we have a set 
of two coupled reaction advection diffusion equations,
\begin{eqnarray}
-\cfrac1{Re} \px^{\,2}
\hUn +u_0\partial_x \hUn +\left\{ s + i\eta v_0 +  \cfrac1{Re}
\eta^2 +  \cfrac1{Re'} \mu_n^2+\cfrac1\eps  B\right\} \hUn  &=&0.
\label{eqprimlinadimn}
\end{eqnarray}
\subsection{Optimal transmission conditions for the baroclinic modes}
The derivation of optimal transmission conditions
for an advection diffusion equation was performed in \cite{martin05}.
Here we are interested in the set of 
coupled reaction advection diffusion equations~\eqref{eqprimlinadimn}.
The baroclinic modes are not coupled with the evolution of the water height.
Hence for these modes the transmission operators have two components
and will be searched on the form
\begin{equation}
\mathcal{B}^n_\pm \ =\  
\left(
\begin{array}{c}
\mp \displaystyle{\frac{1}{Re}} \px \Uh \pm u_0 \Uh 
- \mathcal{S}^{u,n}_\pm\\
\end{array}
\right).
\label{algooptimaln}
\end{equation} 
We search for the solution of system~\eqref{eqprimlinadimn} 
as a sum of exponentials $x\mapsto e^{\lambda x}$. Plugging this ansatz in
the system, we obtain that $ e^{\lambda x}$
solves~\eqref{eqprimlinadimn} if and only if $\lambda$ is a root of
the determinant of the matrix~~ $M_n(\lambda)\ := \ $
\begin{eqnarray*}
\left( 
\begin{array}{cc}
-\cfrac{\lambda^2}{Re}+u_0 \lambda+s+\cfrac{\eta^2}{Re}+\cfrac{\mu_n^2}{Re'}+i\eta v_0&
-\cfrac1\varepsilon\\
&\\
\cfrac1\varepsilon&-\cfrac{\lambda^2}{Re}+u_0 \lambda+s+\cfrac{\eta^2}{Re}+\cfrac{\mu_n^2}{Re'}+i\eta v_0
\end{array}
\right).
\end{eqnarray*}
This determinant is a polynomial of degree four in $\lambda$ and we can compute its four roots
\begin{equation}
\lambda^{n,+}_\pm \ :=\ \frac{Re}{2} \left(u_0 +
\sqrt{\Delta^n_\pm}\right),\ \qquad \ \lambda^{n,-}_\pm \ :=\ \frac{Re}{2} \left(u_0 - \sqrt{\Delta^n_\pm}\right),
\label{lambdan}
\end{equation}
where
\begin{equation}
\Delta^n_\pm \ := \ u_0^2+\frac{4}{Re}\left(\frac{\eta^2}{Re}+\frac{\mu_n^2}{Re'}+s+i\eta v_0\pm\frac{i}{\eps}\right).
\label{delta}
\end{equation}
Every solution $\lambda^{n,\pm}_\pm$ is associated with a one dimensional kernel
generated by the vector $\phi^{n,\pm}_{\pm}$ defined by
\begin{equation*}
\phi^{n,\pm}_{+} = 
\left(
\begin{array}{c}
1\\
-i
\end{array}
\right), \qquad
\phi^{n,\pm}_{-} = 
\left(
\begin{array}{c}
1\\
i
\end{array}
\right).
\end{equation*}
Since the solutions must vanish at infinity 
we search for solutions in $\Omega^-$ on the form
\begin{equation}
\hUnm(x)
\  =\ \alpha^{n,+}_+
  e^{\lambda^{n,+}_+ x} \phi^{n,+}_{+} + \alpha^{n,+}_-
  e^{\lambda^{n,+}_{-} x} \phi^{n,+}_{-}\ =\ 
  \Phi^{n,+} \cdot \exp\left(x \Lambda^{n,+} \right) \cdot \alpha^{n,+},
\label{omegamoinsn}
\end{equation}
where 
\begin{equation}
\Phi^{n,\pm}\ :=\ 
    \left(
    \begin{array}{cc}
1 & 1\\-i&i    
    \end{array}
    \right),
\quad
\Lambda^{n,\pm} \ :=
   \left(
   \begin{array}{cc}
   \lambda^{n,\pm}_+  & 0\\
    0 & \lambda^{n,\pm}_-
    \end{array}
    \right),
\quad
\alpha^{n,\pm}\ :=\  
    \left(
    \begin{array}{c}
     \alpha^{n,\pm}_+\\ 
     \alpha^{n,\pm}_-
    \end{array}
    \right).
\label{Lambda}
\end{equation}
In $\Omega^+$ we search for the solution on the form 
\begin{equation}
\hUnp(x)
\ =\ \alpha^{n,-}_+
  e^{\lambda^{n,-}_+ x} \phi^{n,-}_{+} + \alpha^{n,-}_-
  e^{\lambda^{n,-}_{-} x} \phi^{n,-}_{-}\ = \ 
  \Phi^{n,-} \cdot \exp\left(x\Lambda^{n,-}\right) \cdot \alpha^{n,-}.
\label{omegaplusn}
\end{equation}
It follows from relations~\eqref{omegamoinsn} and~\eqref{omegaplusn} that
\begin{equation}
\px \hUnmp (x) = \Phi^{n,\pm}\cdot  \Lambda^{n,\pm} \cdot \exp\left(x
\Lambda^{n,\pm}\right)\cdot  \alpha^{n,\pm} 
= \Phi^{n,\pm} \cdot \Lambda^{n,\pm} \cdot
\left[\Phi^{n,\pm}\right]^{-1} \cdot \hUnmp.
\label{derivee}
\end{equation}
We can now define the operator $\mathcal{S}^{u,n}_\pm$
in~\eqref{algooptimaln}  in order to derive an optimal algorithm. This
is done through its Laplace-Fourier symbol: 
\begin{equation}
{\hat{\mathcal{S}}}^{u,n}_\pm \ :=\ \mp \frac{1}{Re} \Phi^{n,\pm} \Lambda^{n,\pm} \left[\Phi^{n,\pm}\right]^{-1}
\pm u_0 \mbox{ Id}.
\label{condn}
\end{equation}
\subsection{Approximate transmission conditions for baroclinic modes}
Since we want to construct an efficient but simple Schwarz waveform relaxation algorithm 
we will derive approximated transmission conditions 
by considering an asymptotic analysis of the results of the previous subsection.\\
The definition~\eqref{delta} of $\Delta^n_\pm$ leads to the
expansion ${\hat{\mathcal{S}}}^n_\pm
\ =\ {\hat{\mathcal{S}}}^{u,n}_{\pm,\mbox{ app}}+ O(\sqrt{\eps})$ with
\begin{eqnarray}
{\hat{\mathcal{S}}}^{u,n}_{\pm,\mbox{ app}} 
\ := \frac{1}{2}
\left(
\begin{array}{cc}
\pm u_0 - \displaystyle{\sqrt{\frac{2}{Re}}\frac{1}{\sqrt{\eps}}} 
&  \displaystyle{\sqrt{\frac{2}{Re}}\frac{1}{\sqrt{\eps}}} \\
 & \\
- \displaystyle{\sqrt{\frac{2}{Re}}\frac{1}{\sqrt{\eps}}} & 
\pm u_0 - \displaystyle{\sqrt{\frac{2}{Re}}\frac{1}{\sqrt{\eps}}} 
\end{array}
\right).
\label{snapprox}
\end{eqnarray}
Note that the approximated operator~\eqref{snapprox}
does not depend on $n$. Consequently the related approximated transmission operators~\eqref{algooptimaln}
can be applied to the whole baroclinic velocity, i.e. to the sum of the baroclinic modes.
\subsection{Approximate transmission con{\-}ditions for the baro{\-}tropic mode}
The derivation of optimal transmission conditions
for the linearized viscous shallow water equations 
without advection term was performed in~\cite{martin08}.
Here we are interested in the linearized viscous 
shallow water equations~\eqref{eqprimlinadim0}-\eqref{divlinadim0}.
The transmission operators will be searched on the form~\eqref{algooptimal-}-\eqref{algooptimal+}.\\
\newline
As for the baroclinic modes we search for the solution of system 
\eqref{eqprimlinadim0}-\eqref{divlinadim0}
as a sum of exponentials $e^{\lambda x}$. 
Here $\lambda$ has to be a root of the determinant of the matrix
~$M_0(\lambda)$~ defined by
\begin{eqnarray*}
\left( 
\begin{array}{ccc}
-\cfrac{\lambda^2}{Re}+u_0 \lambda+s+\cfrac{\eta^2}{Re}+i\eta v_0&
-\cfrac1\varepsilon& \lambda/Fr^2\\
&&\\
 \cfrac1\varepsilon&\!\!\!\!-\cfrac{\lambda^2}{Re}+u_0 \lambda+s+\cfrac{\eta^2}{Re}+i\eta v_0 &
\!\!\! \cfrac{i\eta}{Fr^2} \\
&&\\
\lambda& i \eta & s+u_0 \lambda+i\eta v_0 
\end{array}
\right).
\end{eqnarray*}
This determinant is a polynomial of degree five which does not admit a
trivial decomposition. Hence it is not possible to derive an explicit
formula for the solutions
of~\eqref{eqprimlinadim0}-\eqref{divlinadim0}. Consequently, we are not able to obtain an
explicit form for the optimal transmission conditions for the barotropic
mode, even in Fourier-Laplace variables.
In order to derive approximated transmission conditions
we use the fact that the Rossby number is a small parameter 
to compute approximated values of the roots of the determinant of
$M_0(\lambda)$. The related approximated transmission conditions will be coherent with the results of the previous subsection for the baroclinic modes.\\
Since $u_0$ is positive we first notice that three roots~\eqref{lambda00}-\eqref{lambda0-} 
have a negative real part and two roots~\eqref{lambda0+} have a positive real part. 
The {\it negative} roots will be denoted $\lambda^{0,-}_\pm$ and $\lambda^0_0$. 
The {\it positive} ones will be denoted $\lambda^{0,+}_\pm$.
The notations for the related quantities that we introduce later
are coherent with the previous ones~\eqref{Lambda}.
As above, we search 
for the solution in $\Omega^-$ on the form
\begin{eqnarray*}
\hXm(x)
& =& \alpha^{0,+}_+ e^{\lambda^{0,+}_+ x} \phi^{0,+}_{+} 
+ \alpha^{0,+}_- e^{\lambda^{0,+}_{-} x} \phi^{0,+}_{-} 
\ =:\ 
  \Phi^{0,+} \cdot \exp\left(x\Lambda^{0,+}\right) \cdot \alpha^{0,+}.
\end{eqnarray*}
In $\Omega^+$ we search for the solution on the form
\begin{multline*}
\hXp(x)
\ =\ \alpha^{0,-}_+ e^{\lambda^{0,-}_+ x} \phi^{0,-}_{+} 
+ \alpha^{0,-}_- e^{\lambda^{0,-}_{-} x} \phi^{0,-}_{-} 
+ \alpha^{0}_0 e^{\lambda^{0}_{0} x} \phi^{0}_{0} \\
 =:\   \Phi^{0,-} \cdot\exp\left(x\Lambda^{0,-}\right)\cdot \alpha^{0,-}.
\end{multline*}
We compute the following approximations for the roots of the determinant of $M_0(\lambda)$:
\begin{eqnarray}
\lambda^0_0 \, \, \, &=& -\cfrac{s+ i \eta v_0}{u_0} + O(\eps^2),
\label{lambda00}
\\ 
\lambda^{0,-}_\pm&=&  -\cfrac{\sqrt{\pm i Re}}\seps +\left(
\cfrac{Re\,u_0}2 -\cfrac{Re}{4Fr^2u_0} \right) 
+ O(\seps),
\label{lambda0-}
\\
\lambda^{0,+}_\pm&=&  \cfrac{\sqrt{\pm i Re}}\seps+\left(
\cfrac{Re\,u_0}2 -\cfrac{Re}{4Fr^2u_0} \right) 
+ O(\seps).
\label{lambda0+}
\end{eqnarray}
The associated kernel is always one dimensional and spanned by: 
\begin{eqnarray*}
\Phi^0_0 \, \, \, &=& \left(\begin{array}{c}  
0\\0\\1
 \end{array}\right)
+ O(\eps^2),
\\ 
\Phi^{0,-}_\pm&=&\left(\begin{array}{c}  
 u_0 \pm \cfrac{i\sqrt{2Re}}{4Fr^2} \seps\\ \\
\pm i u_0 \mp \cfrac{i\sqrt{2Re}}{4Fr^2} \seps\\ \\
-1 -\left\{ 
\cfrac{ \pm \sqrt2}4 \, \cfrac{i\sqrt{Re}} {u_0Fr^2} 
+\cfrac{\sqrt2}2\, \cfrac{1 \pm i}{\sqrt{Re}} (( \pm u_0+iv_0)\eta+s)
\right\} \seps
 \end{array}\right)
+ O(\eps), 
\\
\Phi^{0,-}_\pm &=&\left(\begin{array}{c}  
 u_0-\cfrac{i\sqrt{2Re}}{4Fr^2} \seps\\ \\
-i u_0 +\cfrac{i\sqrt{2Re}}{4Fr^2} \seps\\ \\
-1 -\left\{ 
\cfrac{-\sqrt2}4\, \cfrac{i\sqrt{Re}} {u_0Fr^2} 
+\cfrac{\sqrt2}2\, \cfrac{1-i}{\sqrt{Re}} ((-u_0+iv_0)\eta+s)
\right\} \seps
 \end{array}\right)
+ O(\eps).
\end{eqnarray*}
As in the baroclinic modes case, we compute the approximated
transmission operators in Laplace-Fourier variables by
\begin{eqnarray*}
{\hat{\mathcal{S}}}^{u,0}_{-,\mbox{ app}}  &=&
\frac{1}{Re} 
\left[
\Phi^{0,-} \Lambda^{0,-} \left[\Phi^{0,-}\right]^{-1}
\right]
_{2,3}
- \left(
\begin{array}{ccc}
u_0 & 0 & -\frac{1}{Re}\\
0 & u_0 & 0
\end{array}
\right),
\\
{\hat{\mathcal{S}}}^{u,0}_{+,\mbox{ app}}  &=&
-\frac{1}{Re} \Phi^{0,+} \Lambda^{0,+} \left[\Phi^{0,+}\right]^{-1} + u_0 Id,
\end{eqnarray*}
where $M_{2,3}$ denotes the first $2 \times 3$ matrix extracted from
the $3 \times 3$ matrix $M$. It leads to the following Laplace-Fourier symbols
\begin{eqnarray}
{\hat{\mathcal{S}}}^{u,0}_{-,\mbox{ app}}=
\frac{1}{2}
\left(\begin{array}{ccccc}
-\cfrac{\sqrt{2}}{\sqrt{Re\eps}}-{u_0} - \cfrac1{Fr^2 u_0} 
&\!&
\cfrac{\sqrt{2}}{\sqrt{Re\eps}}+\cfrac1{2 Fr^2 u_0}
&\!\!&
\cfrac{\!\!\!-2}{Fr^2}
\\ \\
-\cfrac{\sqrt{2}}{\sqrt{Re\eps}}+\cfrac1{2 Fr^2 u_0}
&\!&
-\cfrac{\sqrt{2}}{\sqrt{Re\eps}}-{u_0}
&\!\!&
0
\end{array}\right)
+O(\seps),~
\label{so-approx}
\end{eqnarray}
and
\begin{eqnarray}
{\hat{\mathcal{S}}}^{0}_{+,\mbox{ app}}
=
\frac{1}{2}
\left(\begin{array}{ccc}
-\cfrac{\sqrt{2}}{\sqrt{Re\eps}}+{u_0}-\cfrac1{Fr^2 u_0}
&~&
\cfrac{\sqrt{2}}{\sqrt{Re\eps}}-\cfrac1{2 Fr^2 u_0}
\\ \\
-\cfrac{\sqrt{2}}{\sqrt{Re\eps}}-\cfrac1{2 Fr^2 u_0}
&&
-\cfrac{\sqrt{2}}{\sqrt{Re\eps}}+{u_0}
\\  \\
0&&0
\end{array}\right)
+ O(\seps).
\label{so+approx}
\end{eqnarray}
By using relations~\eqref{snapprox},~\eqref{so-approx} and~\eqref{so+approx}, we notice that
\begin{eqnarray*}
\left[{\hat{\mathcal{S}}}^{u,0}_{\pm,\mbox{ app}}\right]_{2,2}
&=& 
{\hat{\mathcal{S}}}^{u,n}_{\pm,\mbox{ app}}
+
\left(
\begin{array}{ccc}
- \cfrac1{2 Fr^2 u_0} && \mp \cfrac1{4 Fr^2 u_0}\\
\mp \cfrac1{4 Fr^2 u_0} && 0
\end{array}
\right).
\end{eqnarray*}
It follows that a part of the transmission conditions will be applied to the whole velocity
(sum of baroclinic and barotropic modes) while a second part will be applied only to
the barotropic mode. The first part corresponds to the operator ${\hat{\mathcal{S}}}^{u,n}_{\pm,\mbox{ app}}$.
The second one corresponds to the remaining terms in the operator ${\hat{\mathcal{S}}}^{u,0}_{\pm,\mbox{ app}}$.
\subsection{The optimized Schwarz waveform relaxation algorithm}
Thanks to the computed approximated operators~\eqref{snapprox},~\eqref{so-approx} and~\eqref{so+approx}
we can now derive an approximated Schwarz waveform relaxation algorithm for the 
linearized primitive equations~\eqref{eqprimlinadim}--\eqref{inita}.\\
\newline
Since the computed operators~\eqref{snapprox},~\eqref{so-approx} and~\eqref{so+approx}
do not depend neither on the Fourier variable $\eta$ nor on the Laplace variable $s$ 
the related operators in the real space are identical to their Laplace-Fourier symbols.
It follows that the approximated transmission operators $\mathcal{B}_\pm$ 
\eqref{algooptimal-}-\eqref{algooptimal+} 
have the following form
\begin{eqnarray}
\mathcal{B}_- X &=&\left(
\begin{array}{c}
\cfrac1{Re} \partial_x u
+\left(\cfrac{\sqrt{2}}{2\sqrt{Re\,\eps} }-\cfrac{u_0}2\right)u
-\cfrac{\sqrt{2}v}{2\sqrt{Re\,\eps}}+\cfrac{\ub-\vb/2}{2 Fr^2 u_0}\\
\\
\cfrac1{Re} \partial_x v 
+\left(\cfrac{\sqrt{2}}{2\sqrt{Re\,\eps} }-\cfrac{u_0}2\right)v
+\cfrac{\sqrt{2}u}{2\sqrt{Re\,\eps}} -\cfrac{\ub}{4Fr^2 u_0}
\end{array}
\right),
\label{b-}\\
\mathcal{B}_+ X &=&\left(
\begin{array}{c}
-\cfrac1{Re} \partial_x u +\cfrac{\a}{Fr^2}
+\left(\cfrac{\sqrt{2}}{2\sqrt{Re\,\eps} }+\cfrac{u_0}2\right)u
-\cfrac{\sqrt{2}v}{2\sqrt{Re\,\eps}}+\cfrac{\ub+\vb/2}{2 Fr^2 u_0}\\
\\
-\cfrac1{Re} \partial_x v 
+\left(\cfrac{\sqrt{2}}{2\sqrt{Re\,\eps} }+\cfrac{u_0}2\right)v
+\cfrac{\sqrt{2}u}{2\sqrt{Re\,\eps}}+\cfrac{\ub}{4Fr^2 u_0}\\
\\
u_0\a + \ub 
\end{array}
\right),
\label{b+provi}
\end{eqnarray}
for which we recall that $\ub$ and $\vb$ represent the mean-values
with respect to the $z$ variable of the velocities $u$ and $v$.\\
Note that by replacing the first component $({\mathcal{B}_+X})_1$ by
the linear combination  
$$({\mathcal{B}_+X})_1 -1/(Fr^2u_0)({\mathcal{B}_+X)}_3,$$
we replace~\eqref{b+provi} by the equivalent transmission conditions 
\begin{equation}
\mathcal{B}_+^{\sim} X \ =\ \left(
\begin{array}{c}
-\cfrac1{Re} \partial_x u 
+\left(\cfrac{\sqrt{2}}{2\sqrt{Re\,\eps} }+\cfrac{u_0}2\right)u
-\cfrac{\sqrt{2}v}{2\sqrt{Re\,\eps}}-\cfrac{\ub-\vb/2}{2 Fr^2 u_0}\\
\\
-\cfrac1{Re} \partial_x v 
+\left(\cfrac{\sqrt{2}}{2\sqrt{Re\,\eps} }+\cfrac{u_0}2\right)v
+\cfrac{\sqrt{2}u}{2\sqrt{Re\,\eps}}+\cfrac{\ub}{4Fr^2 u_0}\\
\\
u_0\a + \ub 
\end{array}
\right).
\label{b+}
\end{equation}
In the sequel we use~\eqref{b+} rather than~\eqref{b+provi} and we
drop the superscripts ``$\sim$''.\\
Next, we remark that the transmission conditions~\eqref{b-}\eqref{b+} are a particular case of the generalized transmission conditions
\begin{eqnarray}
\mathcal{B}_- X &=&
\cfrac1{Re} \partial_x \Uh
-\cfrac{u_0}2 \Uh
+\cfrac{\alpha}{\sqrt{\eps}} A \Uh + 
\beta B \Ubh,           
\label{b-approxbis}
\\ 
\mathcal{B}_+ X &=&\left(\!\!\!\!
\begin{array}{c}
-\cfrac1{Re} \partial_x \Uh 
+\cfrac{u_0}2 \Uh
+\cfrac{\alpha}{\sqrt{\eps}} A \Uh  
-\beta B \Ubh
\\
u_0\a + \ub 
\end{array}\!\!\!\!
\right)                    
\label{b+approxbis}
\end{eqnarray}
where
\begin{equation}
A:=\left(\begin{array}{cc}1&-1\\1&1\end{array}\right), \qquad
B:=\left(\begin{array}{cc}1&-{1}/{2}\\-{1}/{2}&0\end{array}\right). \label{defAB}
\end{equation}
The original transmission operators~\eqref{b-}-\eqref{b+} correspond to the choice 
\begin{equation}
\label{defab}
\alpha \ =\ \cfrac1{\sqrt{2Re}}\,, \qquad \beta
\ =\   \cfrac{1}{2 Fr^2 u_0}.
\end{equation} 
Notice that $\mathcal{B}_-X$ and $(\mathcal{B}_+X)_{(1,2)}$ do not
depend on the water height $\a$, so we may rewrite
$\mathcal{B}_-X=\mathcal{B}_-^{\Uh} \Uh$ and
$\mathcal{B}_+X=~^t(\mathcal{B}_+^{\Uh} \Uh,\mathcal{B}_+^\a X)$ as 
\begin{equation}
\label{Bpm}
\mathcal{B}_\pm^{\Uh} \Uh\ :=\ \mp \cfrac1{Re} \partial_x \Uh 
\pm\cfrac{u_0}2 \Uh
+\cfrac{\alpha}{\sqrt{\eps}} A \Uh  
\mp \beta B \Ubh 
~\qquad~~~ \mathcal{B}_+^\a X\ :=\ u_0\a + \ub. 
\end{equation}
Let us emphasize the identity: 
\begin{eqnarray}
\label{identite}
\mathcal{B}_+^{\Uh} \Uh+\mathcal{B}_-^{\Uh}\Uh &=&2\cfrac{\alpha}{\sqrt{\eps}} A \Uh.
\end{eqnarray}
This relation will be useful both for defining a weak formulation 
of the algorithm in Section~\ref{secWP} and for the numerical implementation of
this algorithm in Section~\ref{secnum}. \\
Finally the Schwarz waveform relaxation algorithm~\eqref{swra} writes  
\begin{equation}\label{Algo}
\left\{
\begin{array}{rll} \EP(X^{n+1}_-) =&  0&
  \mbox{on }  \Omega^-_T, \\
&&\\
X^{n+1}_-(\cdot,0) =& X_i&\mbox{on } \Omega^- , \\
&&\\
\mathcal{B}_-^{\Uh} \Uhmnn =&\!\!\mathcal{B}_-^{\Uh} \Uhpn&\mbox{on } \Gamma_T, 
\end{array}
          \right.~~~
\left\{
\begin{array}{rll} 
\EP(X^{n+1}_+) =&  0  &\mbox{on } \Omega^+_T,\\
 &&\\
X^{n+1}_+(\cdot,0) =& X_i&\mbox{on } \Omega^+,  \\
&&\\
 \mathcal{B}_+^{\Uh} \Uhpnn =&\!\!\mathcal{B}_+^{\Uh} \Uhmn&\mbox{on }
 \Gamma_T,\\
 &&\\
\mathcal{B}_+^{\a} X_+^{n+1} =&\!\!\mathcal{B}_+^\a X_-^n &\mbox{on }
 \gamma_T.
\end{array}
           \right.
\end{equation}
where the operators $\mathcal{B}_\pm^{\Uh}$, $\mathcal{B}_+^\a$ are defined by
equalities~\eqref{defAB}\eqref{Bpm} and where $\alpha$ and $\beta$ are
free parameters. \\
These generalized transmission conditions can now be optimized with respect to the two parameters $\alpha$ and $\beta$.
In the case of a one dimensional reaction advection diffusion equation
this optimization problem has been solved analytically (see \cite{halpern07}). Here, we will present a numerical procedure in Section \ref{secnum}.
\section{Well-posedness of the linearized Primitive Equations}
\label{secwp}
In the previous sections we have performed formal computations on the linearized
Primitive Equations leading to the construction of the Schwarz
waveform relaxation algorithm~\eqref{Algo}. The aim of this section
is to be more precise: we will define a weak formulation  
of the system~\eqref{eqprimlinadim}--\eqref{inita} and then prove
that this system is well-posed in the natural spaces associated to
this weak formulation. \\
\newline
From now on we relax the boundary condition
on the bottom, i.e. we assume $\alpha_b\geq0$ instead of
$\alpha_b=0$. Moreover, in order to prepare the study of the well posedness of the algorithm~\eqref{swra} in the next section, we
consider non-homogeneous right-hand sides $Y=(F_1,F_2,f)=Y(x,y,z,t)$.
 The system of linearized primitive equations $\EP(X)=Y$ writes
\begin{eqnarray}
\label{un}
\left\{\partial_t  + \U0\cdot \nh -\cfrac1{Re} \Lah  -\cfrac1{Re'}
\pz^{\,2}  +\cfrac1\eps  C\right\}\Uh +\cfrac1{Fr^2} \nh\a&=\ F&\quad\mbox{
  in } \Omega_T,
\\
\label{deux}
\pz \Uh (x,y,0,t) &=\  0&\quad\mbox{ on } \omega_T,\\
\label{trois}
- \pz \Uh (x,y,-1,t) + \alpha_b \Uh (x,y,-1,t) &=\ 0 &\quad\mbox{ on } \omega_T,\\
\label{UN}
\left\{\partial_t  +\U0 \cdot \nh\right\}  \a +\dh \Ubh &=\ f& \quad\mbox{ in } \omega_T.
\end{eqnarray}
We supplement this system with the initial conditions 
\begin{eqnarray}
\label{i}
\Uh(\cdot,0)&=\ \Uhi,& \quad\mbox{ in } \Omega, \\
\label{I}
\a(\cdot,0)&=\ \ai &\quad\mbox{ in } \omega.
\end{eqnarray}
Note that if we consider that the water height $\ap$ is given, 
the system~\eqref{un}--\eqref{trois},~\eqref{i} 
with unknown $\Uh$ is a classical  linear parabolic  problem. On the other hand if we consider that 
the mean horizontal velocities $\Ubh$ are given then $\a$ solves the 
linear transport problem with source term~\eqref{UN},~\eqref{I}.\\
We will proceed as follows: first we recall the classical weak
formulations both for the parabolic problem (with prescribed water height) and
for the transport equations (with prescribed velocity). These two problems
define two maps $\mathcal{S}_1:\ \a \mapsto \Uh$ and $\mathcal{S}_2 :\ \Uh
\mapsto \a$. Finally we define the weak solutions of the Primitive Equations
to be the fixed points of the map $\tau: \ (\Uh,\a)\mapsto
(\mathcal{S}_1(\a),\mathcal{S}_2(\Uh))$ and conclude by proving the existence of
a unique fixed point.\\
\newline
Let us first introduce some functional spaces and some notations. 
We will work with initial data and  right hand sides satisfying 
\begin{eqnarray*}
\Uhi \in H \ := L^2(\Omega ,\R^2), \qquad 
 \ai \in L^2(\omega),\\ 
F\in L^2(0,T;\mathcal{V}'), \qquad f \in L^2(0,T;L^2(\omega)), 
\end{eqnarray*}
where $\mathcal{V}'$ is the topological dual of
$\mathcal{V}:=H^1(\Omega,\R^2)$.\\ 
The weak solutions will satisfy  
\begin{equation*}
\Uh \in \ C\left([0,T], H\right)  \, \cap \,  L^2(0,T;
\mathcal{V}), \quad ~
\a \in \ C\left([0,T],L^2(\omega)\right)\, \cap \,  C(\R_{x},L^2\left(\R_y\times(0,T))\right).
\end{equation*} 
We will need the following bilinear forms:  
\begin{eqnarray}
\nonumber
a(U,V) &:=& \frac{1}{Re} \left( \nabla_h U,  \nabla_h V \right)_{\Omega_T}
+ \frac{1}{Re'} \left(\pz U,\pz V\right)_{\Omega_T} +
\frac{\alpha_b}{Re'} (U,V)_{\omega_{-1,T}}\\
\label{a}
&&\qquad~\qquad  +\frac{1}{\eps}(CU,V)_{\Omega_T}
+(\U0 \cdot \nabla U, V)_{\Omega_T},\\
\label{c}
c(\a,V)&:=&\frac{1}{Fr^2}\left(\a e_x,\px V \right)_{\Omega_T}.
\end{eqnarray}
where $\omega_{-1}:=\R_x\times \R_y \times \{-1\}_z$ and
$(U,V)_\Sigma$ denotes the $L^2$ scalar product on $\Sigma$.\\ 
\newline
Assuming that we have a strong solution, and taking the scalar product
of equation~\eqref{un} with $V\in
\mathcal{D}(\overline{\Omega}\times(0,T),\R^2)$,  we obtain (after integrating by parts) the following weak formulation
for the equations governing the horizontal velocities:
\begin{equation}
\label{wun} 
\forall\,V\in \mathcal{D}(\overline{\Omega}\times(0,T),\R^2),\qquad
(\partial_t \Uh,V)_{\Omega_T}^{}+a(\Uh,V)=c(\a ,V)+\lra{ F,V}.
\end{equation}
We now state\\

\bd
Let $F\in L^2(0,T;\mathcal{V}')$ and $\a\in L^2(\omega_T)$, we say that $\Uh\in
L^2(0,T;\mathcal{V})$ is a weak solution of the system~\eqref{un} if~\eqref{wun} holds. 
\ed

\medskip
 
\bp
\label{wpun}
Let $\Uhi\in L^2(\Omega)$, $F\in L^2(0,T;\mathcal{V}')$ and $\a\in
L^2(\omega_T)$, there exists a unique weak solution $\Uh\in C([0,T];H)\cap
L^2(0,T;\mathcal{V})$ of~\eqref{un} satisfying the initial
condition~\eqref{i}. Moreover, we have the energy inequality 
\begin{multline}
\label{energie-un}
\frac{1}{2} \|\Uh\|^2_{\Omega}(t) + \int_0^t \left\{
\cfrac1{Re}\|\nabla_h \Uh\|^2_{\Omega}(s) +\cfrac1{Re'}\|\pz
\Uh\|^2_{\Omega}(s) + \cfrac{\alpha_b}{Re'} \|\Uh\|_{\omega_{-1}}^2
(s)\right\}\,ds \\ 
\leq \  \frac{1}{2} \|\Uhi\|^2_{\Omega}+ \int_0^t\left\{ \lra{ F\,,\,\Uh
}(s) + (\partial_x \ub,\a)_\omega(s) \right\}\, ds. \end{multline}
\ep

\medskip

\begin{proof}
The method is classical and we only sketch the proof.
We obtain the existence of a solution satisfying~\eqref{energie-un} by
the Galerkin method.
Let $(E_m)$ be an increasing sequence of
\textit{finite dimensional} sub-spaces of  
$\mathcal{V}$ such that $\cup E_m$ is dense in $\mathcal{V}$. For
every $m$, there exists $U_m\in C^\infty([0,T],E_m)$ such that~\eqref{wun} holds for
every $V\in \mathcal{D}(0,T;E_n)$ -- we only have to solve a finite
system of linear ordinary differential equations.\\
Using $U_m\times\bs{1}_{[0,t]}$ as a test function in~\eqref{wun}, 
we conclude that $U_m$ satisfies~\eqref{energie-un}. 
So using the Cauchy Schwarz inequality and the Gr\"onwall Lemma, 
we see that the sequence $(U_m)$ is uniformly bounded in $L^2(0,T;\mathcal{V})$. \\
Now from the weak formulation, we deduce that $(\partial_t U_m)$ is bounded in
$L^2(0,T;\mathcal{V}')$, thus by Aubin-Lions Lemma, $(U_m)$ is
compact in $C([0,T],H)$. Extracting a subsequence we obtain a
solution satisfying~\eqref{energie-un}.\\
Regularizing in time and using the weak formulation, we see that any solution
satisfies~\eqref{energie-un} and uniqueness follows by the energy
method. 
\end{proof}

\medskip

\noindent
Let us turn our attention to the equations~\eqref{UN},~\eqref{I} governing the
evolution of the water height $\a$. This is a linear transport
equation with constant coefficients and a source
term. Assuming that $\a$ is a strong solution, multiplying~\eqref{UN}
by a test function $\b\in \mathcal{D}(\omega\times [0,T))$, integrating
on $\omega_T$, integrating by parts in space \textit{and} time and then using the initial
condition~\eqref{I}, we obtain 
\begin{multline}
\label{wUN}
 \forall \b \in \mathcal{D}(\omega \times[0,T)),\\
-\left(\a ,  \left\{\partial_t+ \U0\cdot \nabla_h\right\} \b
\right)_{\omega_T}=\ \left(\ai,\b(\cdot,0) \right)_{\omega} + (f-\nabla_h\cdot\Ubh,\b)_{\omega_T}. 
\end {multline}

\medskip

\bd
Let $f\in L^2(\omega_t)$, $\Uh \in L^2(0,T;\mathcal{V})$ and $\ai\in
L^2(\omega)$, we say that $\a\in L^2(\omega_T) $ is a weak solution of the
system~\eqref{UN},~\eqref{I} if~\eqref{wUN} holds. 
\ed

\medskip

\noindent
Remark that the test function does not necessarily vanish at time $0$
and that the initial data is prescribed by the weak formulation.

\medskip

\bp
\label{wpUN}
Let $f\in L^2(\omega_t)$, $\Uh \in
L^2(0,T;\mathcal{V})$ and $\ai\in L^2(\Omega)$.  There exists a unique
weak  solution $\a \in L^2(\omega_T)$
of~\eqref{UN},~\eqref{I}. Moreover this solution is given
by the characteristic formula:  
\begin{equation} 
\label{solUN}
\a(x,y,t)=\ai(x-u_0t,y-v_0t)+\int_0^t(f- \nabla_h\cdot
\Ubh)(x-u_0 s,y-v_0 s,t-s)ds.
\end{equation}
This solution lies in 
$C\left([0,T];L^2(\omega)\right)\, \cap\, C\left(\R_x;L^2(\R_y \times
(0,T))\right)$  and satisfies the following estimates for every $t \in
[0,T]$ and every $x\in \R$, 
\begin{eqnarray}
\|\a(\cdot,\cdot,t)\|_{\omega}&\leq &\|\ai\|_{\omega} 
+ \int_0^t \|f-\nabla_h \cdot\Ubh\|_{\omega}(s)ds,
\label{estim1}\\
\|\a(x,\cdot,\cdot)\|_{\gamma_t}&\leq& \frac{1}{u_0}
\left(\|\ai\|_{\omega}+ \int_0^t \|f-\nabla_h\cdot\Ubh\|_{\omega}(s)ds\right).
\label{estim2}
\end{eqnarray}
\ep

\medskip

\begin{proof}
First, notice that the  estimates \eqref{estim1}~\eqref{estim2} are direct
consequences of the characteristic formula~\eqref{solUN}.\\
Next, remark that if the data $\nabla_h \Ubh, f$ and $\ai$
are sufficiently smooth then the function $\a$ given by the
formula~\eqref{solUN} solves~\eqref{wUN}. Hence we obtain the existence of a
solution of~\eqref{wUN} by density. 
\\
For the uniqueness, by linearity we may assume that the data $\nabla_h
\Ubh, f$ and $\ai$ vanish. Then let $\psi \in \mathcal{D}(\omega)$ and
$\rho \in \mathcal{D}([0,T))$ and define the test function $\b$ by
$\b(x,y,t)\ :=\ -\psi(x-u_0 t,y-v_0 t) \int_t^T \rho(s)ds$, so that:
$$
\partial_t \b + \U0 \cdot \nabla  \b = \psi(x-u_0t,y-v_0 t)\rho(t),$$ 
and~\eqref{wUN} yields 
\begin{equation*}
0=\int_{\omega_T} \a(x,t)\rho(t)\psi(x-u_0 t,y-v_0 t)
\ =\ \int_{\omega_T} \a(x+u_0 t,y+v_0 t,t)\rho(t)\psi(x,y).
\end{equation*}
Since this is true for every  $(\psi,\rho) \in
\mathcal{D}(\omega)\times\mathcal{D}([0,T))$, we have $\a\equiv 0$ on $\omega_T$.
\end{proof}

\medskip

\noindent
Finally, we define the notion of weak solution for the linearized
primitive equations.  

\medskip

\bd
Let $Y=(F,f)\in L^2(\Omega_T,\mathcal
{V}')\times L^2(\omega_T)$ and
$X_i=(\Uhi,\ai)\in L^2(\Omega)\times L^2(\omega)$. We say that $X=(\Uh,\a)\in C(0,T;H)\times L^2(\omega_T)$ is a weak
solution of~\eqref{un}--\eqref{I} if the weak
formulations~\eqref{wun} and~\eqref{wUN} hold and if $\Uh(\cdot,0)=\Uhi$.
\ed

\medskip

\bt
\label{thm1}
Let $Y=(F,f)\in L^2(\Omega_T,\mathcal
{V}')\times L^2(\omega_T)$ and
$X_i=(\Uhi,\ai)\in L^2(\Omega)\times L^2(\omega)$. There exists a
unique weak solution $X=(\Uh,\a)\in (C(0,T;H)\cap L^2(0,T;\mathcal{V}))\times L^2(\omega_T)$ of~\eqref{un}--\eqref{I}.
\et

\medskip

\begin{proof}
The right hand side $Y$ and the initial data $X_i$ being fixed,
Proposition~\ref{wpun} and Proposition~\ref{wpUN} define two maps
\begin{eqnarray*}
  S_1:\ L^2(\omega_T)&\rightarrow& C(0,T;H)\cap
L^2(0,T;\mathcal{V}),\qquad \a \ \mapsto\  \Uh,
\end{eqnarray*}
and
\begin{eqnarray*}
 S_2:\ L^2(0,T,\mathcal{V})&\rightarrow&
C\left([0,T];L^2(\omega)\right) \cap C\left(\R_x;L^2(\R \times
(0,T))\right),\qquad \Uh \ \mapsto\ \a.
\end{eqnarray*}
Denoting by $\mathcal{T}$ the
affine mapping  $(\Uh,\a) \mapsto (S_1(\a),S_2(\Uh))$, the application
$X$ is a weak solution of~\eqref{un}--\eqref{I} if and only if it is a fixed point of $\mathcal{T}$ in $$ \mathcal{E}_T:= L^2(0,T;\mathcal{V}) \times C([0,T],L^2(\omega)).$$  
Let $X_1,X_2 \in \mathcal{E}_T$ and let $(\Uh,\a):=X_1-X_2$ and
$(\pUh,\pa):=\mathcal{T}(X_1)-\mathcal{T}(X_2)$, by linearity,
using~\eqref{energie-un},
 we get for $0\leq t\leq T$,
\begin{multline*}
\frac{1}{2} \|\pUh\|^2_{\Omega}(t) + 
\int_0^t \left\{\cfrac1{Re} \|\nabla_h \pUh\|^2_{\Omega}(s) +\cfrac1{Re'}\|\pz
\pUh\|^2_{\Omega}(s)\right\}\,ds \leq \int_0^t (\partial_x
\overline{\tilde{u}},\a)_\omega(s) ds\\
 \leq \ \left( \int_0^t \| \nabla \pUh\|_{\Omega}^2(s)
 \,ds\right)^{1/2} \left(\int_0^t \|\a\|^2_{\omega}(s)\, ds\right)^{1/2} .
\end{multline*}
By Young inequality, we may absorb the term in $\nabla \pUh$ in the
left hand side and get:
\begin{equation}
\label{preuve_wp1}
\|\pUh\|^2_{\Omega}(t) +  \int_0^t\|\nabla \pUh\|_{\Omega}^2(s)\,ds \ \leq\ \kappa t \sup_{s\in[0,t]} \{\|\a\|^2_{\omega}(s)\} \qquad \mbox{for } 0\leq t\leq T,
\end{equation}
for some $\kappa >0$. Now~\eqref{estim1} and the Cauchy Schwarz
inequality yield 
\begin{equation}
\label{preuve_wp2}
\|\pa\|^2_{\omega}(t)\ \leq \ t \int_0^t \|\nabla 
 \Uh\|_{\Omega}^2(s)ds, \qquad \mbox{for } 0\leq t\leq T.
\end{equation}
Finally, inequalities~\eqref{preuve_wp1}~\eqref{preuve_wp2} imply that, for $T'\in (0,T]$  small enough, the mapping $\mathcal{T}$ is strictly contracting in
$\mathcal{E}_{T'}$ yielding the existence of a unique fixed point of
$\mathcal{T}$ in $\mathcal{E}_{T'}$. Repeating the argument on the
intervals $[T',2T']$, $[2T',3T']$, ... we obtain the result on $[0,T]$.
\end{proof}

\section{Weak formulation and well-posedness of the Schwarz waveform relaxation algorithm}
\label{secWP}
We study in this section the well-posedness of the
algorithm~\eqref{Algo}. First, we will define weak
formulations for the two sub-problems and prove that they are well-posed.  We will pay a particular attention to the weak
form of the transmission conditions. In particular we will establish
that the solutions $X^{n+1}_\pm$ of the $n^{th}$ step  of the 
algorithm~\eqref{Algo} are in the right spaces, allowing the
construction of the transmission conditions for the next step.\\
\newline  
As in the previous section, we also consider non-homogeneous
right-hand sides $Y=(F,f)$. Every step of the algorithm 
may be split  in the two following sub-problems.  First in the domain $\{x<0\}$, we search
for a solution $X^{n+1}_-:=X_-=(\Uhm,\am)$ solving the initial and boundary
value parabolic problem, 
\begin{equation}
\label{Algo-un}
\left\{.
\begin{array}{rl}
\left\{\partial_t  +\U0\! \!\cdot\! \nh  -\cfrac1{Re} \Lah-\cfrac1{Re'}
\pz^{\,2}  +\cfrac1\eps  C\right\}\Uhm +\cfrac1{Fr^2} \nh\am\ =\ F
&\mbox{ in }\Omega^-_T,\\
\\
 - \pz \Uhm (x,y,-1,t) + \alpha_b \Uhm (x,y,-1,t) \ =\ 0&\!\!\!\!,  \\
\pz \Uhm (x,y,0,t) \ =\ 0 &\mbox{ on } \omega_T^-, \\
\\
\mathcal{B}_-^{\Uh} \Uhm\ =\ \mathcal{B}_-^{\Uh} \Uhpn&\mbox{ on }
\Gamma_T, \\
\\
\Uhmnn(\cdot,0) \ =\ \Uhi &\mbox{ in } \Omega^- , \\
\end{array}\right.
\end{equation}
and the transport problem, 
\begin{equation}
\label{Algo-UN}
\left\{
\begin{array}{rcl}
\left\{\partial_t  +\U0 \cdot \nh\right\} \am +\dh \Ubhm =&f \quad&\mbox{ in }\omega^-_T,\\
\amnn(\cdot,0) =&\ai &\mbox{ in } \omega^-.
\end{array}\right.
\end{equation}
In the right subdomain $\{x>0\}$ we search for a solution $X^{n+1}_+:=X_+=(\Uhp,\ap)$  solving the initial and boundary
value parabolic problem, 
\begin{equation}
\label{Algo+un}
\left\{.
\begin{array}{rl}
\left\{\partial_t  +\U0\! \!\cdot\! \nh  -\cfrac1{Re} \Lah-\cfrac1{Re'}
\pz^{\,2}  +\cfrac1\eps  C\right\}\Uhp +\cfrac1{Fr^2} \nh\ap\ =\ F
&\mbox{ in }\Omega^+_T,\\
\\
 - \pz \Uhp (x,y,-1,t) + \alpha_b \Uhp (x,y,-1,t) \ =\ 0&\!\!\!\!,  \\
\pz \Uhp (x,y,0,t) \ =\ 0 &\mbox{ on } \omega_T^+, \\
\\
\mathcal{B}_+^{\Uh} \Uhp\ =\ \mathcal{B}_+^{\Uh} \Uhmn&\mbox{ on }
\Gamma_T, \\
\\
\Uhp(\cdot,0) \ =\ \Uhi &\mbox{ in } \Omega^+ , \\
\end{array}\right.
\end{equation}
and the transport problem with entering characteristics on the boundary $\gamma_T$,
\begin{equation}
\label{Algo+UN}
\left\{
\begin{array}{rcl}
\left\{\partial_t +\U0 \cdot \nh \right\} \am +\dh \Ubhp =&f \quad&\mbox{ in }\omega^+_T,\\
\mathcal{B}_+^{\a} X_+\ =&\mathcal{B}_+^\a X_-^n &\mbox{ on } \gamma_T,\\
\ap(\cdot,0)\ =&\ai &\mbox{ in } \omega^+.
\end{array}\right.
\end{equation}
To prove that these two sub-problems are well-posed, we proceed as in Section~\ref{secwp}.
First we study the  parabolic problems with prescribed water heights:
we introduce a weak formulation for these problems and prove that they
are well-posed. Then we study the transport equations, introduce their
weak formulations and establish their well-posedness. Finally, the
solutions of the coupled parabolic-transport problems are obtained \textit{via}
a fixed point method.\\
\newline
As in Section~\ref{secwp}, the initial data $X_i(\Uhi,\ai)$ satisfy
$\Uhi\in H$, $\ai \in L^2(\omega)$. We choose right hand sides $Y=(F,f)$ in $L^2(0,T;H)\times
L^2(\omega_t)$. (In section~\ref{secwp}, we only assumed $F\in L^2(0,T;\mathcal{V}')$, but here this choice would cause
difficulties at the interface). We will search for weak solutions $X_\pm=(\Uhpm,\apm)$ in
the spaces,    
\begin{eqnarray}
\label{regularite1}
\Uhpm &\in &C\left([0,T], H^\pm\right)  \quad\cap\quad  L^2(0,T;
\mathcal{V}^\pm),\\
\label{regularite2}
\apm &\in &C\left([0,T],L^2(\omega^\pm)\right)\quad \cap\quad  C(\R_{\pm,x},L^2\left(\R_y\times(0,T)_t)\right),
\end{eqnarray} 
with $H^\pm\ :=\ L^2(\Omega^\pm,\R^2)$ and  $\mathcal{V}^\pm\ :=\ H^1(\Omega^\pm,\R^2)$.
\subsection{The parabolic problems}
Let us define the weak-formulation for the parabolic
problems~\eqref{Algo-un} and~\eqref{Algo+un}. First we introduce the bilinear forms $a^\pm$ and $c^\pm$:  
\begin{eqnarray}
\nonumber
a^\pm(U,V) &:=& \frac{1}{Re} \left( \nabla_h U,  \nabla_h V \right)_{\Omega^\pm_T}
+ \frac{1}{Re'} \left(\pz U,\pz V\right)_{\Omega^\pm_T} +
\frac{\alpha_b}{Re'} (U,V)_{\omega^\pm_{-1,T}}\\
\label{apm}
&&\qquad~\qquad  +\frac{1}{\eps}(CU,V)_{\Omega^\pm_T}
+(\U0 \cdot \nabla U, V)_{\Omega^\pm_T},\\
\label{cpm}
c^{\pm}(\a,V)&=&\frac{1}{Fr^2}\left(\a e_x,\px \overline{V} \right)_{\omega^\pm_T} 
  \pm \frac{1}{Fr^2} \left(  \a e_x,  \overline{V}\right)_{\gamma_T},
\end{eqnarray}
where $\omega^\pm_{-1}:=\R_x^\pm \times \R_y \times \{-1\}_z$.\\
Next, taking the scalar product of the first equation of~\eqref{Algo-un}
or~\eqref{Algo+un} with some test map $V\in
\mathcal{D}(\overline{\Omega^\pm}\times(0,T),\R^2)$,  we obtain:
\begin{equation*} 
(\partial_t \Uhpm,V)_{\Omega^\pm\times(0,T)}^{}+a^\pm(\Uhpm,V)=c^\pm(\apm ,V)
  \mp \frac{1}{Re}(\px \Uhpm , V)_\Gamma^{}+ (F,V)_{\Omega^\pm_T}.
\end{equation*}
Then, using the transmission conditions to express $\partial_x \Uhpm$ on $\Gamma$, we get
\begin{multline*} 
(\partial_t \Uhpm,V)_{\Omega^\pm\times(0,T)}^{}+a^\pm(\Uhpm,V)+b^\pm(\Uhpm,V)\\
=  c^\pm(\apm ,V) + \left(\mathcal{B}^{\Uh}_\pm\Uhmpn ,V
\right)_{\Gamma}+ (F,V)_{\Omega^\pm_T}.
\end{multline*}
with  
\begin{equation}
\label{bpm}
b^\pm (U,V) \ :=\ \pm \cfrac{u_0}2 (U , V)_\Gamma^{}+ \cfrac{\alpha}{\sqrt{\eps}} (A U ,  V)_\Gamma^{} \mp \beta (B \overline{U} ,  \overline{V})_\Gamma^{}.
\end{equation}\\
We are still not satisfied with this weak formulation. Indeed, the
knowledge of $\px \Uhmpn$ on the boundary $\Gamma\times(0,T)$  is
needed for defining the term $(\mathcal{B}^{\Uh}_\pm\Uhmpn
,V)_{\Gamma}$ in the right hand side. 
Unfortunately,~\eqref{regularite1} only gives: $\px \Uhmpn \in
L^2(\Omega^\pm\times(0,T))$ which is not sufficient to define a trace. To overcome this
difficulty, we use relation~\eqref{identite} to define recursively the
terms $(\mathcal{B}^{\Uh}_\pm\Uhmpn ,V)_{\Gamma}$. Indeed, for strong
solutions, we have on $\Gamma_T$
\begin{equation*}
\mathcal{B}^{\Uh}_\mp\Uhpm  \ \stackrel{\eqref{identite}}{=}\  
  -\mathcal{B}^{\Uh}_\pm\Uhpm  +2\cfrac{\alpha}{\sqrt{\eps}}  A \Uhpm
  \ \stackrel{\eqref{Algo}}{=}\  -\mathcal{B}^{\Uh}_\pm \Uhmpn
  +2\cfrac{\alpha}{\sqrt{\eps}}  A \Uhpm.
\end{equation*}
Thus, identifying $\mathcal{B}^{\Uh}_\pm\Uhmpn$ with a distribution  
$\mathcal{B}^n_\pm\in L^2(0,T;\mathcal{W}')$, where $\mathcal{W}$ denotes the space $H^{1/2}(\Gamma,\R^2)$ ; we
obtain a weak formulation of the algorithm for the horizontal velocities:

\medskip

\bd
Assuming that the functions $\apm=\apmnn$ are known, the weak formulation
of the parabolic part~\eqref{Algo-un} and~\eqref{Algo+un} of the
algorithm~\eqref{Algo}  are defined as follows:\\
For the first step, we choose
\begin{equation}
\label{AlgoW1}
\mathcal{B}^0_\pm\in L^2(0,T;\mathcal{W}') 
\end{equation}
Then for $n\geq 0$, the horizontal velocity is defined by $\Uhpmnn=\Uhpm$ where
$\Uhpm$ solves 
\begin{multline} 
\forall\,V\in \mathcal{D}(\overline{\Omega^\pm}\times(0,T),\R^2),\quad
(\partial_t \Uhpm,V)_{\Omega^\pm\times(0,T)}^{}+a^\pm(\Uhpm,V)+b^\pm(\Uhpm,V)\\
=  c^\pm(\apm ,V) + \left< \mathcal{B}^n_\pm\,,\,V\right>_{\Gamma_T} + (F,V)_{\Omega_T^\pm},
\label{AlgoW2}
\end{multline}
where $a^\pm$, $c^\pm$, and $b^\pm$ are defined in~\eqref{apm}---\eqref{bpm}.
Once $\Uhpmnn$ is known, we can define the boundary conditions
for the next step in the opposite domain by
\begin{equation}
\label{AlgoW3}
\mathcal{B}_\mp^{n+1} \ :=\  - \mathcal{B}_\pm^{n}  + 2 \cfrac{\alpha}{\sqrt{\eps}} A {\Uhpmnn}_{|\Gamma_T}.
\end{equation}
\ed 
Notice that assuming that the maps $\Uhpmnn$
satisfy~\eqref{regularite1} then their traces on $\Gamma_T$ are well
defined in $L^2(0,T;\mathcal{W})\subset L^2(0,T;\mathcal{W}')$. Consequently, 
the transmission conditions $\mathcal{B}_\mp^{n+1}$ defined
recursively by~\eqref{AlgoW3} stay in the space $L^2(0,T;\mathcal{W}')$.\\

\bp
\label{wp_Algo_pm}
Let $\Uhi\in H$, $F\in L^2(0,T;H)$, $\mathcal{B}_\pm^n\in
L^2(0,T;\mathcal{W}_\pm')$ and $\apm$ $(=\apmnn)$ satisfying~\eqref{regularite2}. Then there exists a
unique $\Uhpmnn=\Uhpm$  with regularity~\eqref{regularite1}
satisfying~\eqref{AlgoW2} and the initial condition
$\Uhpmnn(0)\equiv\Uhi$ on $\Omega_\pm$. Moreover, we have the energy inequality 
\begin{multline}
\label{Energie+-un}
\frac{1}{2} \|\Uhpm\|^2_{\Omega^\pm}(t) +\left(\cfrac{\alpha}{\sqrt{\eps}} \pm {u_0}/2\right) \|\Uhpm\|_{\Gamma_t}^{2} 
\mp \beta \left(\|\ub_\pm\|_{\gamma_t}^2-(\ub_\pm,\vb_\pm)_{\gamma_t}\right)\\
+ \int_0^t \left\{
\cfrac1{Re}\|\nabla_h \Uhpm\|^2_{\Omega^\pm}(s) +\cfrac1{Re'}\|\pz
\Uhpm\|^2_{\Omega^\pm}(s) + \cfrac{\alpha_b}{Re'} \|\Uhpm\|_{\omega^\pm_{-1}}^2
(s)\right\}\,ds\\
\leq \  \frac{1}{2} \|\Uhi\|^2_{\Omega^\pm}+
(F\,,\,\Uhpm)_{\Omega_t^\pm}+ \lra{  \mathcal{B}_\pm^{n} ,\Uhpm }_{\Gamma_t} \\
+\int_0^t \left\{(\partial_x \ub_\pm,\apm)_{\omega^\pm}(s) \pm ( \ub_\pm,\apm)_{\gamma^\pm}(s)\right\}\, ds.
\end{multline}
\ep
\begin{proof}
We proceed as in the proof of Proposition~\ref{wpun}: we apply the
Galerkin method.  Here we only check that the~\textit{a priori} 
inequality~\eqref{Energie+-un} is sufficient for applying this
method.  In order to bound the quadratic terms in the left
hand side of~\eqref{Energie+-un} and the last term in the right hand side, we will use the inequality 
\begin{eqnarray*}
\|U\|_{\Gamma}^2 &\leq & 2 \|U\|_{\Omega^\pm} \| \px U\|_{\Omega_\pm},
 \end{eqnarray*}
valid for   $U \in \mathcal{V}_\pm$. (To prove it, write 
$|U(0,y,z)|^2= 2\int_{-\infty}^0(\partial_x U\cdot U)(x',y,z) \, dx'$ 
integrate on $\R_y\times(-1,0)_z$ and use the Cauchy-Schwarz inequality). 
From this inequality, the  Cauchy-Schwarz inequality, the Young
inequality and the fact that the trace on $\Gamma$ defines a continuous embedding
$\Pi\ :\ \mathcal{V}^\pm \rightarrow \mathcal{W}$, we see that~\eqref{Energie+-un} implies 
\begin{multline}
\label{Energie+-deux}
\|\Uhpm\|^2_{\Omega^\pm}(t) + \int_0^t \|\nabla
\Uhpm\|_{\Omega^\pm}^2(s)\, ds -\kappa \int_0^t \|\Uhpm\|^2_{\Omega^\pm}(s)\, ds\\ 
\leq \ \kappa \left\{ \|\Uhi\|_{\Omega^\pm}^2 +
\|F\|_{\Omega^\pm_t}^2 + \int_0^t \|
\mathcal{B}_\pm^{n}\|_{\mathcal{W}'}^2(s)\, ds +
\|\apm\|^2_{\Omega_t} + \|\apm\|^2_{\gamma_t}\right\}
\end{multline}
for some $\kappa>0$. Taking a Galerkin sequence $(U_m)$ associated
to~\eqref{AlgoW2}, the elements of this sequence
satisfy~\eqref{Energie+-un} and then inequality~\eqref{Energie+-deux} and the Gr\"onwall
Lemma imply that this sequence is bounded in
$L^2(0,T;\mathcal{V})$. Extracting a subsequence (as in Proposition~\ref{wpun}) we obtain a
solution of~\eqref{AlgoW2}. \\
Then using the weak formulation satisfied by $U_m$ we see that $(\partial_t U_m)$ is bounded in $L^2(0,T;\mathcal{V}')$ and from Aubin-Lions Lemma (see e.g.~\cite{SW97}), the sequence $(U_m)$ is compact $L^2(0,T; H^s(\Omega^\pm,\R^2))$ for
$s<1$. Thus we may let $m$ tend to $\infty$ in the quadratic boundary terms in the left hand side  of~\eqref{Energie+-un}. \\
The uniqueness follows from~\eqref{Energie+-deux} and Gr\"onwall Lemma.
\end{proof}
\subsection{The transport equations}
We now consider that the velocities $\Uhpm=\Uhpmnn$ are known and study the
transport problems~\eqref{Algo-UN} and~\eqref{Algo+UN}. We begin with
the domain $\{x<0\}$. Proceeding exactly as in Section~\ref{secwp}, we
 obtain that a strong solution of Problem~\eqref{Algo-UN} satisfies 
\begin{multline}
\label{w-UN}
 \forall \b \in \mathcal{D}(\omega^- \times[0,T)),\\
-\left(\am ,  \left\{\partial_t+ \U0\cdot \nabla_h\right\} \b
\right)_{\omega^-_T}=\ \left(\ai,\b(\cdot,0) \right)_{\omega^-} + (f-\nabla_h\cdot\Ubhm,\b)_{\omega^-_T}. 
\end {multline}
\bd
Let $f\in L^2(\omega_t)$, $\Uhm(=\Uhmnn) \in L^2(0,T;\mathcal{V}^-)$ and $\ai\in
L^2(\omega)$. We say that $\am\in L^2(\omega_T^-) $ is a weak solution
of Problem~\eqref{Algo-UN} if~\eqref{w-UN} holds. 
\ed

\medskip

The following result is proved exactly as Proposition~\ref{wpUN}\\

\bp
\label{PropTransport}
Let $f\in L^2(\omega_t)$, $\Uhmnn \in L^2(0,T;\mathcal{V}^-)$ and $\ai\in
L^2(\omega)$.  There exists a unique
weak  solution $\amnn =\am\in L^2(\omega_T^-)$
of~\eqref{Algo-UN}. Moreover this solution is explicitly given by the formula:  
\begin{equation} 
\am(x,y,t)=\ai(x-u_0t,y-v_0t)+\int_0^t(f- \nabla_h\cdot \Ubhm)(x-u_0 s,y-v_0 s,t-s)ds.
\end{equation}
It lies in 
$C\left([0,T];L^2(\omega^-)\right)\, \cap\, C\left((-\infty,0]_x;L^2(\R_y \times
(0,T))\right)$  and satisfies the following estimates for every $t \in
[0,T]$ and every $x\leq 0$, 
\begin{eqnarray}
  \|\am(\cdot,t)\|_{\omega^-}\leq &\|\ai\|_{\omega^-} 
+ \int_0^t \|f-\nabla_h \cdot\Ubhm\|_{\omega^-}(s)ds,&
\label{estim-1}
\\
\|\am(x,\cdot)\|_{\gamma_t}\leq& \frac{1}{u_0}
\left(\|\ai\|_{\omega^-}+ \int_0^t \|f-\nabla_h\cdot\Ubhm\|_{\omega^-}(s)ds\right).&
\label{estim-2}
\end{eqnarray}
\ep
Once the solutions of~\eqref{Algo-un}-\eqref{Algo-UN} are known it is
possible to define the transmission conditions on the water-height for
the next step  (see~\eqref{Bpm}) 
\begin{equation}
 \label{AlgoW4}
\mathcal{B}_+^\a X_-^{n+1} \ :=\ u_0 \amnn(0,\cdot))+ \ubmnn(0,\cdot).
\end{equation}
In the domain ${x>0}$, the situation is slightly different since there are
ingoing characteristics on $\gamma_T$. So we choose test functions that
do not necessarily vanish on the boundary and use the transmission
condition to prescribe the value of the solution on
$\gamma_T$. Finally, a solution of~\eqref{Algo+UN} satisfies 
\begin{multline}
\label{w+UN}
 \forall \b \in \mathcal{D}(\overline{\omega}^+ \times[0,T)),\quad 
-\left(\ap ,  \left\{\partial_t+ \U0\cdot \nabla_h\right\} \b \right)_{\omega^+_T}\\=\ \left(\ai,\b(\cdot,0) \right)_{\omega^+} +\left(\ab,\b(0,\cdot) \right)_{\R_t}   + (f-\nabla_h\cdot\Ubhp,\b)_{\omega^+_T}, 
\end {multline}
where the boundary value $\ab$ is defined on $\gamma_T$ by
\begin{equation}
\label{w+DEUX}
\ab \ :=\ \cfrac1{u_0} \left\{\mathcal{B}_+^\a X_-^n -\ub_+ \right\}.
\end{equation}
\bd
\label{defi101}
Let $f\in L^2(\omega_t)$, $\Uhp$ $(=\Uhpnn) \in L^2(0,T;\mathcal{V}^+)$, $\ai\in
L^2(\omega)$. Assuming that $\ab$ defined by~\eqref{w+DEUX} belongs to
$L^2(\gamma_T)$ , we say that $\apnn=\ap\in L^2(\omega_T^+) $ is a weak
solution of Problem~\eqref{Algo+UN} if~\eqref{w+UN} holds. 
\ed

\medskip

\noindent
Using the characteristic method, we have\\

\bp
Let $f$, $\Uhp$ $(=\Uhpnn)$, $\ai$ and $\ab$ be as in Definition~\ref{defi101}. There exists a unique
weak  solution $\apnn \in L^2(\omega_T^+)$
of~\eqref{Algo+UN}. Moreover it is given by the
characteristic formula:  
\begin{equation*} 
\ap(x,y,t)\ =\ \ai(x-u_0t,y-v_0t)+\int_0^t(f- \nabla_h\cdot \Ubhp)(x-u_0
s,y-v_0s,t-s)ds
\end{equation*}
if  $x> u_0 t$, and
\begin{equation*} 
\ap(x,y,t)\ =\ \ab\left(y-\frac{v_0}{u_0}x ,t-\frac{x}{u_0}\right)+\int_0^{\frac{x}{u_0}}(f- \nabla_h\cdot \Ubhp)(x-u_0
s,y-v_0s,t-s)ds,
\end{equation*}
with $\ab$ given by \eqref{w+DEUX}, if $x\leq  u_0 t$. \\
\newline
The solution belongs to 
$C\left([0,T];L^2(\omega^+)\right)\, \cap\, C\left([0,+\infty)_x;L^2(\R_y \times
(0,T))\right)$  and satisfies the following estimates for every $t \in
  [0,T]$ and every $x\geq 0$, 
\begin{eqnarray}
  \|\ap(\cdot,t)\|_{\omega^+}\leq &\|\ai\|_{\omega^+} +  u_0 \|\ab\|_{\gamma_t} 
+ \int_0^t \|f-\nabla_h \cdot\Ubhp\|_{\omega^+}(s)ds,&
\label{estim+1}
\\
\|\ap(x,\cdot)\|_{\gamma_t}\leq& \frac{1}{u_0}
\left(\|\ai\|_{\omega^+}+  u_0\|\ab\|_{\gamma_t} +\int_0^t \|f-\nabla_h\cdot\Ubhp\|_{\omega^+}(s)ds\right).&
\label{estim+2}
\end{eqnarray}
\ep
\subsection{Well-posedness of the algorithm}
First we define a weak formulation for the left and right sub-problems
at step $n$ of the algorithm.\\

\bd
\label{defsubprob}
Let $Y=(F,f)\in L^2(\Omega_T)\times L^2(\omega_T)$, let
$X_i=(\Uhi,\ai)\in L^2(\Omega)\times L^2(\omega)$. 
For $n\geq 0$. 
\begin{itemize}
\item Let $\mathcal{B}_-^n\in L^2(0,T;\mathcal{W}')$. Then
$X^{n+1}_-=(\Uhm,\am)$ is a weak solution of
Problem~\eqref{Algo-un},~\eqref{Algo-UN} if it has
regularity~\eqref{regularite1}-\eqref{regularite2} and if $\Uhm$ 
(respectively $\am$)  is  a weak solution of~\eqref{Algo-un}(respectively~\eqref{Algo-UN}). 
\item Let $\mathcal{B}_+^n\in L^2(0,T;\mathcal{W}')$ and
  $\mathcal{B}_+^\a X^n_- \in  L^2(\gamma_T)$. Then  $X^{n+1}_+=(\Uhp,\ap)$ is a weak solution of
Problem~\eqref{Algo+un},~\eqref{Algo+UN} if it has
regularity~\eqref{regularite1}-\eqref{regularite2} and if $\Uhp$,
(respectively $\ap$)  is  a weak solution of~\eqref{Algo+un}(respectively~\eqref{Algo+UN}). 
\end{itemize} 
\ed

\medskip
 
\noindent
Then we give a weak formulation for the complete algorithm.\\

\bd
\label{WeakAlgo}
The weak formulation of Algorithm~\eqref{Algo} is defined by  
\begin{itemize}
\item Choose $\mathcal{B}_+^\a X^0_- \in L^2(\gamma_T)$ and
  $\mathcal{B}_\pm^0\in L^2(0,T;\mathcal{W}')$.  
\end{itemize} 
Then, for $n\geq 0$,
\begin{itemize}
\item Find $X_-^{n+1}$ weak
  solution of~\eqref{Algo-un}-\eqref{Algo-UN} 
and $X_+^{n+1}$ weak
  solution of~\eqref{Algo+un}-\eqref{Algo+UN}. 
\item Define the transmission conditions for step $n+1$ by~\eqref{AlgoW3} and~\eqref{AlgoW4}.  
\end{itemize}
\ed

\medskip

\bt
\label{thm2}
With the hypotheses of Definition~\ref{defsubprob}, there exists a
unique weak solution $X_-^{n+1}$ (respectively $X_+^{n+1}$)
of~Problem~\eqref{Algo-un},\eqref{Algo-UN} (respectively~\eqref{Algo+un},\eqref{Algo+UN}).  
\et

\medskip

\begin{proof}
We only prove the result for the left sub-problem, the other one being
similar. As in the proof of Theorem~\ref{thm1}, we use a fixed point
method. Let us introduce the spaces 
\begin{eqnarray*}
\mathcal{E}^1_T& :=& C\left([0,T], H^-\right)  \quad\cap\quad  L^2(0,T;
\mathcal{V}^-),\\
\mathcal{E}^2_T&:=& C\left([0,T],L^2(\omega^-)\right)\quad \cap\quad  C((-\infty,0]_{x},L^2\left(\R_y\times(0,T)_t)\right).\end{eqnarray*}
Proposition~\ref{wp_Algo_pm} (respectively
Proposition~\ref{PropTransport}) defines an affine mapping
$\mathcal{S}_1:\ \mathcal{E}^2_T \rightarrow \mathcal{E}^1_T $, $\am
\mapsto \Uhm$ (respectively
$\mathcal{S}_2:\  \mathcal{E}^1_T\rightarrow \mathcal{E}^2_T$, 
$\Uhm\mapsto\am$).\\
Setting $\mathcal{E}^-_T:=\mathcal{E}^1_T\times \mathcal{E}^2_T-$, an
application $X_-^{n+1}$ is a weak solution of
Problem~\eqref{Algo-un},\eqref{Algo-UN} if and only if it is a fixed
point in $\mathcal{E}^-_T$ of the mapping 
$$
\mathcal{T}^-\ :\ (\Uhm,\am)\mapsto(\mathcal{S}_1^-(\am),\mathcal{S}_2^-(\Uhm)).
$$
We now show that $\mathcal{T}^-$ has a unique fixed point. Let $X_1,X_2
\in \mathcal{E}_T^-$ and let $(\Uhm,\am):=X_1-X_2$ and
$(\pUh,\pa):=\mathcal{T}^-(X_1)-\mathcal{T}^-(X_2)$. By
linearity~\eqref{Energie+-deux} yields: for $0\leq t\leq T$,
\begin{equation*}
\|\pUhm\|^2_{\Omega^-}(t) + \int_0^t \|\nabla
\pUhm\|_{\Omega^-}^2(s)\, ds -\kappa \int_0^t \|\pUhm\|^2_{\Omega^-}(s)\, ds 
\ \leq \ \kappa \left\{
\|\am\|^2_{\Omega_t} + \|\am\|^2_{\gamma_t}\right\}.
\end{equation*}
And from Gr\"onwall lemma, we obtain for $ 0\leq t\leq T$,
\begin{multline}
\label{preuve-wp1}
\|\pUhm\|^2_{\Omega}(t) +  \int_0^t\|\nabla \pUhm\|_{\Omega}^2(s)\,ds
\\ \leq\ \kappa e^{\kappa T}\left\{ t \sup_{[0,t]} \|\am(\cdot,s)\|^2_{\omega} + \sup_{\R_-}  \|\am(w,\cdot)\|^2_{\gamma_t}  \right\}.
\end{multline}
Now from~\eqref{estim-1} and~\eqref{estim-2}, we get 
\begin{equation}
\label{preuve-wp2}
\|\pam\|^2_{\omega^-}(t)+u_0 \|\pam (x,\cdot)\|_{\gamma_t}\ \leq \ t \|\nabla 
 \Uhm\|_{\Omega^-_t}^2 \qquad \mbox{for } 0\leq t\leq T.
\end{equation}
Finally, we endow $\mathcal{E}^-_t$ with the norm $\|(\Uhm,\am)\|_{\mathcal{E}^-_t}:=$ 
\begin{equation*}\left( \sup_{[0,t]}
\|\Uhm(\cdot,s)\|^2_{\Omega^-} +\|\nabla \Uhm\|^2_{\Omega_t^-} +
\sup_{[0,t]} \|\am(\cdot,s)\|^2_{\omega}
+ 2\kappa e^{\kappa T} \sup_{\R_-} \|\am(x,\cdot)\|^2_{\R_t}\right)^{1/2}.
\end{equation*}
With this norm~\eqref{preuve-wp1}~\eqref{preuve-wp2} imply that for
$T'\in (0,T]$  small enough, $\mathcal{T}^-$ is contracting in
 $\mathcal{E}_{T'}$.  This yields the existence of a unique fixed point of
$\mathcal{T}^-$ in $\mathcal{E}_{T'}^-$. We obtain the result on
$[0,T]$ by continuation.
\end{proof}

\medskip

\noindent
Finally, we can state\\

\bt
\label{thm3}
The algorithm~\eqref{WeakAlgo} is well-defined.
\et

\medskip

\begin{proof}
We only have to check that for each step the hypotheses of
Theorem~\ref{thm2} are satisfied. The solutions $X_\pm^{n+1}$, build
at step $n$ have
regularity~\eqref{regularite1},~\eqref{regularite2}. We easily deduce
that  $\mathcal{B}_\pm^{n+1}$ defined by~\eqref{AlgoW3} belongs to
$L^2(0,T;\mathcal{W}')$ and $\mathcal{B}_+^\a X_-^{n+1}$ defined by~\eqref{AlgoW4}
belongs to $L^2(\gamma_T)$. Thus the hypotheses of Theorem~\ref{thm2}
hold for step $n+1$.   
\end{proof}

\medskip

\br
\label{rem_convergence}
Although we do not exhibit a proof here, we are able to establish the convergence of the algorithm in some cases. More precisely, if the matrices $A$ and $B$ defined by~\eqref{defAB} are replaced by diagonal matrices $\tilde{A}$ and $\tilde{B}$, $\tilde{A}$ being positive definite and $\tilde{B}$ being non negative, then the algorithm converges. The proof relies on the energy method developed for the Shallow water equations without advection term in~\cite{martin08}. 
Modifying slightly the proof, we can allow $\tilde{A}$ and $\tilde{B}$ to have non vanishing skew-symmetric off-diagonal parts. This generalization still does not cover the situation~\eqref{defAB} because $B$ has a symmetric non vanishing off-diagonal part. 
Nevertheless, numerical evidences of the convergence of the algorithm are given in the next section.  
\er

\section{Numerical results}
\label{secnum}
\subsection{Numerical scheme in the subdomains}
\begin{center}
\begin{figure}[hbtp]
\rotatebox{0}{\includegraphics[scale=.65]{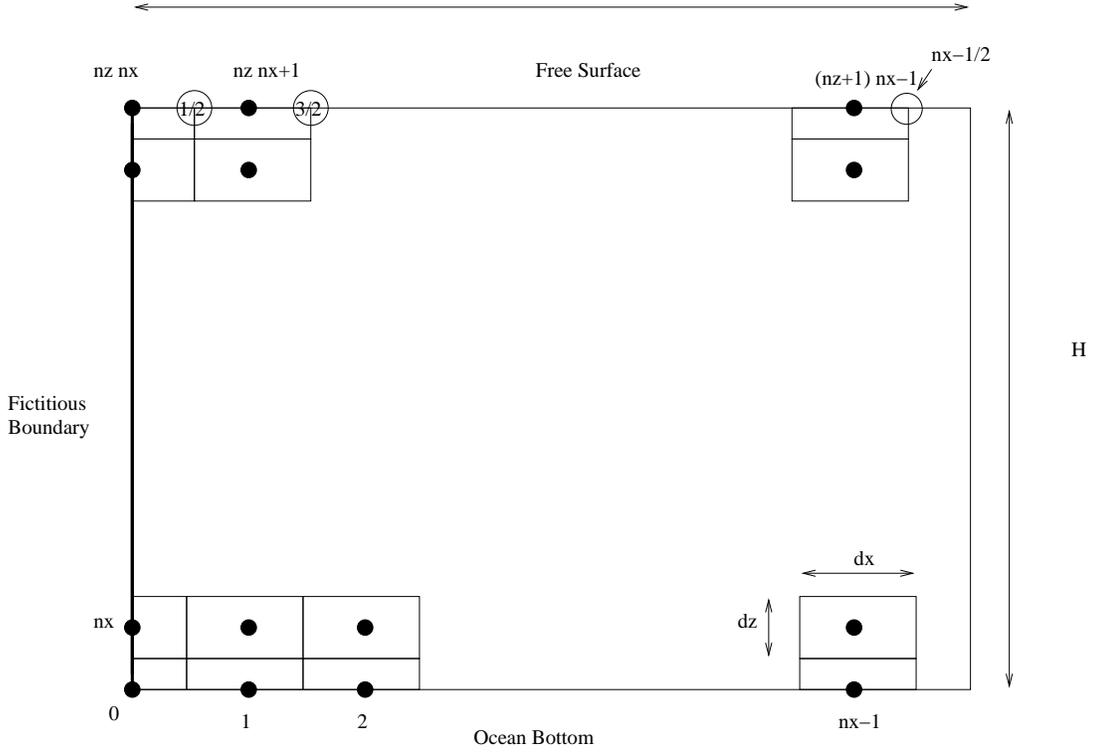}}
\caption{\label{fig:omega+}Space discretization of $\Omega^+$}
\end{figure}
\end{center}

For the numerical applications we consider for simplicity a 2 dimensional domain
and the related two dimensional $(x,z)$ version 
of the primitive equations (\ref{eqprimlinadim})-(\ref{inita}).
Note that the transmission conditions (\ref{b-approxbis})-(\ref{b+approxbis}) 
are independent of the transverse $y$-variable 
and are not affected by this simplification.\\
\newline
In this subsection we do not deal with the boundary conditions.
Hence the processes are the same in both subdomains $\Omega^\pm$
and we restrict ourselves to the subdomain $\Omega^+$.
We first describe the space discretization of the subdomains. 
We consider a regular cartesian grid 
of $nx \times nz$ points and we apply a finite volume method.
We introduce the horizontal space step $\Delta x$ 
and the vertical space step $\Delta z$.
For Euler or Navier-Stokes type problems 
it is well known that a good way to recover some numerical stability 
is to compute velocities and pressure on different cells
(see for instance \cite{arakawa} and the publications devoted to the so-called C-grids).
Here we only deal with the horizontal velocity and the water height (depending  only depending on $x$ and $t$) plays the role of the pressure. 
We thus have to introduce two types of finite volume meshes - see Figure \ref{fig:omega+}.
The first one is a 2d finite volume mesh
and is related to the computation of the velocities.
For $i=0... n_x-1$ and $j=0... n_z$ we denote $I=i+j n_x$.
The cells of this first mesh will be denoted 
$C^+_{I} = X^+_{I} + (-\Delta x / 2, \Delta x / 2) \times (-\Delta z / 2, \Delta z / 2)$.
where the points $X^+_{I}$ stand for  $X^+_{I} = (0,-H) + (i\Delta x,j\Delta z)$ 
(they are represented by a black circle in Figure \ref{fig:omega+}).
The second grid is a 1d finite volume mesh devoted to the computation of the water height.
The cells of this second mesh will be denoted 
$c_{i+1/2} = x_{i+1/2} + (-\Delta x / 2, \Delta x / 2)$.
where the points 
$x_{i+1/2}$ stand for $x_{i+1/2} = (i+1/2) \Delta x$
(they are represented by a circle with a number inside in Figure \ref{fig:omega+}).
Let us now consider the discretization of the equations.
Let us start with momentum equation (\ref{eqprimlinadim}). 
We integrate it on the time-space cell $[t_k,t_{k+1}] \times C_{I}$. 
We compute the interface fluxes at time $t_{k+1/2}$ by classical centered formulas.
We recover the well-known Crank-Nicolson scheme.
It is known to be second order accurate 
and conditionally stable in the $L^\infty$ norm 
under a CFL type condition on the time step $\Delta t_k = t_{k+1}-t_k$.
This strategy is applied for all the velocity nodes such that
the neighboring nodes are included inside the considered subdomain.
The related discrete relations stand
\begin{multline}
{u_{I,k+1}} 
+ \frac{\dt}{2} \bigg\{u_0 \Dx u_{I,k+1} - \cfrac1{Re} 
\Dx^2 u_{I,k+1}  
-\cfrac1{Re'} \Dz^2 u_{I,k+1} 
- \cfrac1\eps v_{I,k+1}  + \cfrac1{Fr^2 \dx}
\overline{\Dx}_1\,\a_{j,k+1}\bigg\}\\ =
{u_{I,k}} 
 - \frac{\dt}{2} \bigg\{u_0 \Dx u_{I,k} - \cfrac1{Re} 
\Dx^2 u_{I,k}  
-\cfrac1{Re'} \Dz^2 u_{I,k} 
- \cfrac1\eps v_{I,k}  + \cfrac1{Fr^2}
\overline{\Dx}_1\,\a_{j,k}\bigg\}.
\end{multline}
\begin{multline}
{v_{I,k+1}} 
 + \cfrac\dt2 \bigg\{u_0 \Dx v_{I,k+1} - \cfrac1{Re} 
\Dx^2 v_{I,k+1}  
-\cfrac1{Re'} \Dz^2 v_{I,k+1} 
+ \cfrac1\eps u_{I,k+1}  \bigg\}\\
=\ {v_{I,k}} - \cfrac\dt2 \bigg\{u_0 \Dx v_{I,k} - \cfrac1{Re} 
\Dx^2 v_{I,k}  
-\cfrac1{Re'} \Dz^2 v_{I,k} 
+ \cfrac1\eps u_{I,k} \bigg\}.
\end{multline}
where $\Dx u_{I,k} = (u_{I+1,k}-u_{I-1,k})/\dx$ 
denotes a classical approximation of the first derivative in space in horizontal direction,
$\Dx^2 u_{I,k} = (u_{I+1,k}-2u_{I,k}+u_{I-1,k})/(2\dx)$ 
and $\Dz^2 u_{I,k} = (u_{I+nx,k}-2u_{I,k}+u_{I-nx,k})/(2\dx)$ 
denote classical approximations of second derivatives in space
in horizontal and vertical directions, respectively.\\
\newline
Let us now consider the mass equation (\ref{divlinadim}).
We integrate it on time space cells $[t_k,t_{k+1}] \times c_{i+1/2}$
- except for $i=0$ where we need to use the transmission conditions.
We compute the interface fluxes by using explicit upwind formulas.
The resulting scheme is known to be first order 
and also conditionally stable under a CFL type condition.
The related formula stands
\begin{eqnarray*}
\a_{i+1/2,k+1} 
&=&  \left(1-\cfrac{\dt}\dx \, u_0 \right) \a_{i+1/2,k}
+\cfrac{\dt}\dx \, u_0 \, \a_{i-1/2,k} \\
&& \quad- \cfrac{\dt\dz}\dx \left(
\ub_{i,k}^{+,n+1}-{\ub_{i-1,k}^{+,n+1}}
\right)
\nonumber
\end{eqnarray*}
where 
$$
\ub_{i,k}^{+,n+1}= 
\cfrac{u_{j,k}}2
+ \sum_{j=1}^{n_z-1} u_{jn_x+j,k} + \cfrac{u_{n_zn_x+j,k}}2
$$
denotes the discrete mean velocity of the flow along the vertical direction. 
\subsection{Numerical discretization near the boundaries}
We now have to explain how we compute the numerical solution 
when one of the interfaces of the cell belongs to the physical boundaries
of the domain or to the fictitious one that is related to the domain decomposition method.
For all cases we choose to work in the same finite volume framework that we use in the interior of the subdomains.\\
\newline
For the physical boundary conditions (\ref{boundlinadim})  
we use the ghost cells method.
This method consists in introducing a fictitious cell along the boundary
and then using the same finite volume strategy as in the interior domain.
For no slip conditions (\ref{boundlinadim}) 
we choose the values 
of the unknowns in the fictitious cell to be equal 
to their values in the neighboring interior cell.\\
\newline
Let us now focus on the numerical treatment 
on the cells that are connected with the interface $\Gamma= \partial \Omega^+ \cap \partial \Omega^-$ -
see Fig. \ref{fig:omega+}. Here we will use a discrete version of the transmission conditions (\ref{b-approxbis})-(\ref{b+approxbis}). This discrete information 
will be the only data that will be transmitted from a subdomain to the other one.
Let us first consider the mass equation (\ref{divlinadim}).
We integrate it on the cell $c_{3/2}$ to obtain
$$
\dx\left[\a^{+,n+1}_{1/2,k+1} -\a^{+,n+1}_{-1/2,k}\right] + \dt\left(
u_0\left[\a^{+,n+1}_{1/2,k} - {\a_{-1/2,k}^{+,n+1}}  \right]  + \ub_{1,k}^{+,n+1}-{\ub_{0,k}^{+,n+1}}\right) = 0.
$$
where $\a^{+,n+1}_{1/2,k+1}$ denotes the water height computed in cell $c^+_{1/2}$ at time $t^{k+1}$
and for iteration $n+1$ of the algorithm. The quantities ${\a_{-1/2,k}^{+,n+1}}$ and ${\ub_{0,k}^{+,n+1}}$ 
have to be considered as unknown quantities since the corresponding cells are not included in $\Omega^+$.
We will use the transmission conditions (\ref{Algo}) to evaluate them.
Hence we obtain thanks to (\ref{Bpm})
$$
\a_{1/2,k+1}^{+,n+1} =\left(1-\cfrac{u_0\dt}\dx\right)\a_{1/2,k}^{+,n+1}
- \cfrac\dt\dx \ub_{1,k}^{+,n+1}
+\cfrac\dt\dx {\mathcal{B}}^{\a,n}_{+,k}
$$
where ${\mathcal{B}}^{\a,n}_{+,k}$ has been computed in $\Omega^-$
during the previous Schwarz iteration and is given by
$$
{\mathcal{B}}^{\a,n}_{+,k}=u_0\a^{-,n}_{n_x,k}
+\ub_{n_x,k}^{-,n}
$$
The basic idea is the same for the momentum equation (\ref{eqprimlinadim}).
Here we integrate the equation on the semi-cell 
${\tilde C}^+_{I} = X^+_{I} + (0, \Delta x / 2) \times (-\Delta z / 2, \Delta z / 2)$.
for $I=jn_x$ with $j=0,...,n_z$. We obtain
\begin{eqnarray*}
&&\cfrac{\dx\dz}2 (u_{I,k+1}^{+,n+1} - u_{I,k}^{+,n+1})\\
&&+ \dt\dz \Bigg\{\cfrac{u_0}2
\left(\cfrac{u_{I,k}^{+,n+1}(r)+u_{I,k+1}^{+,n+1}(r)}2 - {u_{I,k+1/2}^{+,n+1}(l)}\right)\\
&&\qquad \qquad -\cfrac1{Re}\left(\cfrac{\Dx u_{I,k}^{+,n+1}(r)+\Dx
  u_{I,k+1}^{+,n+1}(r)}2-{\partial_xu_{I,k+1/2}^{+,n+1}(l)}\right)\\
&&\qquad \qquad +\cfrac1{Fr^2}\left(\cfrac{\a_{1/2,k}^{+,n+1}+\a_{1/2,k+1}^{+,n+1}}2- {\a^{+,n+1}_{-1/2,k+1/2}}\right)\Bigg\} \\
&&\qquad \qquad-\cfrac{\dt\dx\dz}2\ \cfrac1{Re'}\ \cfrac{ \Dz^2 u_{I,k}^{+,n+1} + \Dz^2 u_{I,k+1}^{+,n+1}}2
-\cfrac{\dt\dx\dz}2 \ \cfrac1\eps \ \cfrac{v_{I,k}^{+,n+1}+v_{I,k+1}^{+,n+1}}2 \ =\ 0
\end{eqnarray*}
where $l$ (respectively $r$) denotes quantities that are evaluated 
on the left (respectively right) boundary of the cell ${\tilde C}^+_{I}$. 
Hence quantities ${u_{I,k+1/2}^{+,n+1}(l)}$, ${\partial_xu_{I,k+1/2}^{+,n+1}(l)}$
and ${\a^{+,n+1}_{-1/2,k+1/2}}$ are unknown quantities since they involve quantities 
that are computed outside the domain $\Omega^+$. Here also we use transmission conditions 
(\ref{Algo}) 
and we obtain thanks to (\ref{Bpm})
\begin{eqnarray}
\label{b21}
&&\cfrac{u_{I,k+1}^{+,n+1}}\dt 
 + \cfrac12 \bigg\{u_0 \cfrac{u_{I,k+1}^{+,n+1}(r)}{\dx} - \cfrac1{Re} 
\cfrac{2\Dx u_{I,k+1}^{+,n+1} (r)}{\dx}  -\cfrac1{Re'} \Dz^2 u_{I,k+1}^{+,n+1} 
- \cfrac1\eps v_{I,k+1}^{+,n+1} \bigg\}\nonumber\\
&&\qquad- \cfrac1\dx\left\{-\cfrac{\alpha}{\sqrt{\eps}} u^{+,n+1}_{I,k+1} + \cfrac{\alpha}{\sqrt{\eps}} v^{+,n+1}_{I,k+1}
-(\cfrac1{Fr^2u_0}-\beta) \bar{u}^{+,n+1}_{1,k+1} - \cfrac\beta2 \bar{v}^{+,n+1}_{1,k+1}\right\} \nonumber\\
&&\qquad= -\cfrac2\dx
{\mathcal{B}}^{u,n}_{+,j,k}
+\cfrac1\dx\left\{-\cfrac{\alpha}{\sqrt{\eps}} u^{+,n+1}_{I,k} +\cfrac{\alpha}{\sqrt{\eps}} v^{+,n+1}_{I,k}
-(\cfrac1{Fr^2u_0}-\beta) \bar{u}^{+,n+1}_{j,k} - \cfrac\beta2 \bar{v}^{+,n+1}_{j,k}\right\}
\nonumber\\
&&\qquad \qquad \qquad
 - \cfrac12 \bigg\{u_0 \cfrac{u_{I,k}^{+,n+1}(r)}{\dx} 
 - \cfrac1{Re} \ \cfrac{2\Dx u_{I,k}^{+,n+1}(r)}{\dx} 
-\cfrac1{Re'} \Dz^2 u_{I,k}^{+,n+1} 
- \cfrac1\eps v_{I,k}^{+,n+1} \bigg\}\nonumber\\
&&\qquad \qquad \qquad +\cfrac{u_{I,k}^{+,n+1}}\dt  -\cfrac1{Fr^2} \cfrac{\a^{+,n+1}_{1/2,k}+\a^{+,n+1}_{1/2,k+1}}\dx  
\end{eqnarray}
where ${\mathcal{B}}^{u,n}_{+,j,k}$ has been computed in $\Omega^-$
during the previous Schwarz iteration and is deduced from relation (\ref{identite})
$$
{\mathcal{B}}^{u,n}_{+,j,k}={\mathcal{B}}^{u,n-1}_{-,j,k} 
+ 2 \cfrac{\alpha}{\sqrt{\eps}} \cfrac{u^{-,n}_{j n_x,k}+u^{-,n}_{j n_x,k+1}}{2} 
- 2 \cfrac{\alpha}{\sqrt{\eps}} \cfrac{v^{-,n}_{j n_x,k}+v^{-,n}_{j n_x,k+1}}{2}
$$
Note that ${\mathcal{B}}^{u,n-1}_{-,j,k}$ is known since it has been computed in $\Omega^+$
at iteration $n-1$ and has been transmitted to the domain $\Omega^-$ before iteration $n$.
Same type of computations for the transverse component of the velocity lead to the following scheme
\begin{eqnarray}
\label{b22}
&&\cfrac{v_{I,k+1}^{+,n+1}}\dt 
 + \cfrac12 \bigg\{u_0 \cfrac{v_{I,k+1}^{+,n+1}(\dx)}{\dx} 
 - \cfrac1{Re} \cfrac{2\Dx^+ v_{I,k+1}^{+,n+1}}{\dx}  
-\cfrac1{Re'} \Dz^2 v_{I,k+1}^{+,n+1} 
+ \cfrac1\eps u_{I,k+1}^{2,n+1} \bigg\}\\
&&\qquad- \cfrac1\dx\left\{-\cfrac{\alpha}{\sqrt{\eps}} u^{+,n+1}_{I,k+1} -\cfrac{\alpha}{\sqrt{\eps}} v^{+,n+1}_{I,k+1}
-\cfrac\beta2 \bar{u}^{+,n+1}_{1,k+1}\right\} \nonumber\\
&&\qquad= -\cfrac2\dx
{\mathcal{B}}^{v,n}_{+,j,k}+\cfrac1\dx\left\{-\cfrac{\alpha}{\sqrt{\eps}} u^{+,n+1}_{I,k} -\cfrac{\alpha}{\sqrt{\eps}} v^{+,n+1}_{I,k}
-\cfrac\beta2 \bar{u}^{2,n+1}_{j,k}\right\}\nonumber\\
&&\qquad \qquad+\cfrac{v_{I,k}^{+,n+1}}\dt 
 - \cfrac12 \bigg\{u_0 \cfrac{v_{I,k}^{+,n+1}(r)}{\dx} 
 - \cfrac1{Re} \cfrac{2\Dx v_{I,k}^{+,n+1}(r)}{\dx} -\cfrac1{Re'} \Dz^2 v_{I,k}^{+,n+1} 
+ \cfrac1\eps u_{I,k}^{+,n+1} \bigg\}\nonumber
\end{eqnarray}
where ${\mathcal{B}}^{v,n}_{+,j,k}$ has been computed in $\Omega^-$
during the previous Schwarz iteration and is deduced from relation (\ref{identite})
$$
{\mathcal{B}}^{v,n}_{+,j,k}= {\mathcal{B}}^{v,n-1}_{-,j,k}  
+2\cfrac{\alpha}{\sqrt{\eps}} \cfrac{u^{-,n}_{j n_x,k}+u^{-,n}_{j n_x,k+1}}{2} 
+2\cfrac{\alpha}{\sqrt{\eps}} \cfrac{v^{-,n}_{j n_x,k}+v^{-,n}_{j n_x,k+1}}{2} 
$$
The derivation of the discrete boundary condition in $\Omega^-$
is based on the same type of computations. 
Note that only the components of the velocity 
are concerned by the transmission problem in $\Omega^-$.
\subsection{Numerical optimization of the transmission conditions}
\begin{figure}[hbtp]
\centering
\includegraphics[width=.5\textwidth,angle=-90]{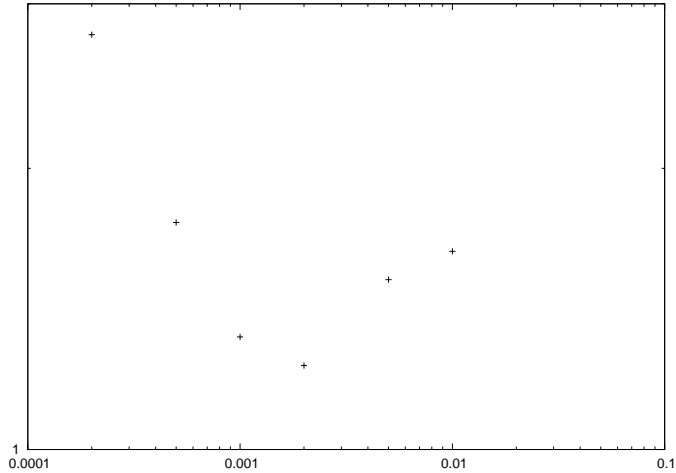}
\caption{\label{fig:alphaepsilon} Quotient $\alpha_{opt}/\alpha_{Tay}$ between the numerically optimized parameter and the Taylor approximation parameter as a function of the Rossby number $\eps$ (in Log scale)}
\end{figure}

\begin{figure}[hbtp]
\centering
\includegraphics[height=.35\hsize,angle=0]{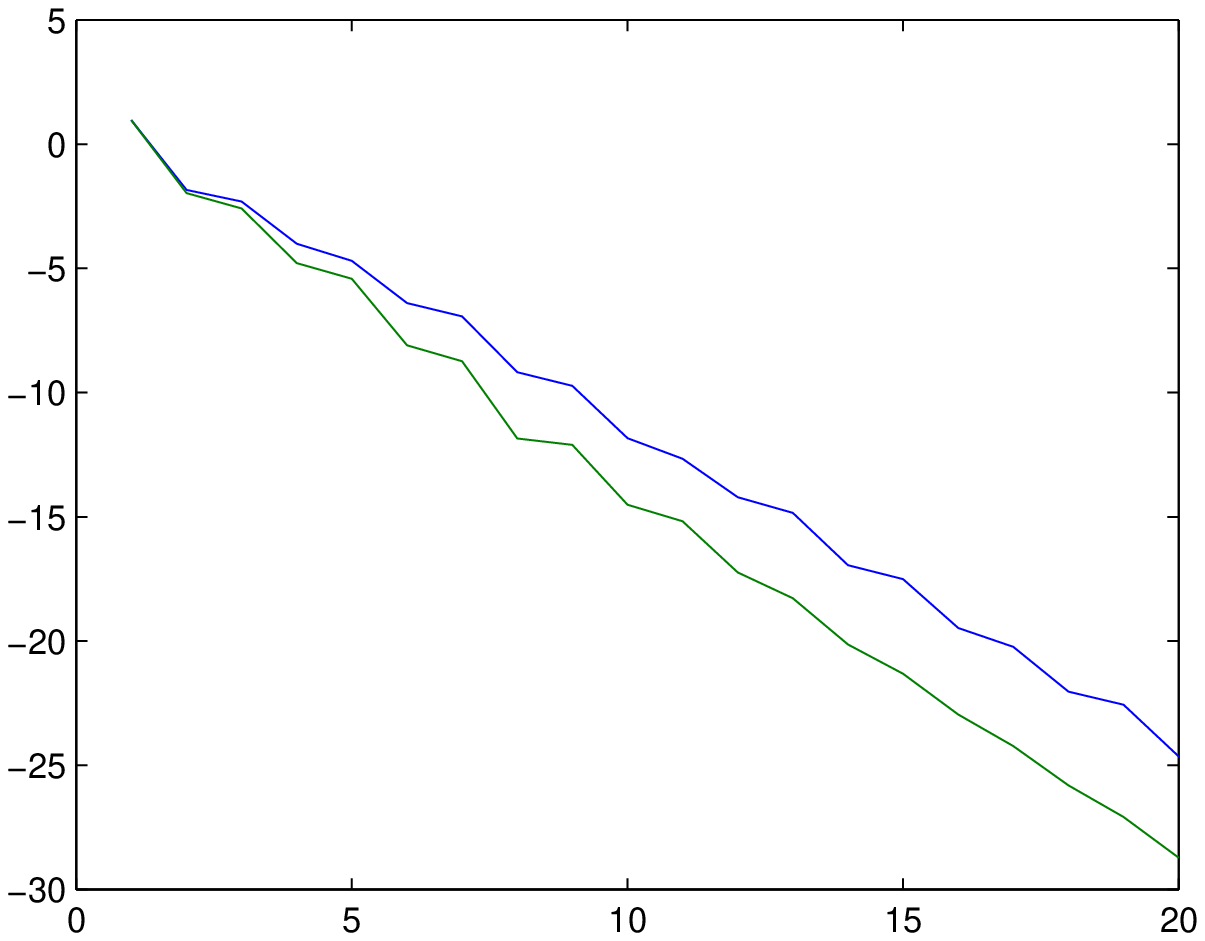}
\ 
\includegraphics[height=.35\hsize,angle=0]{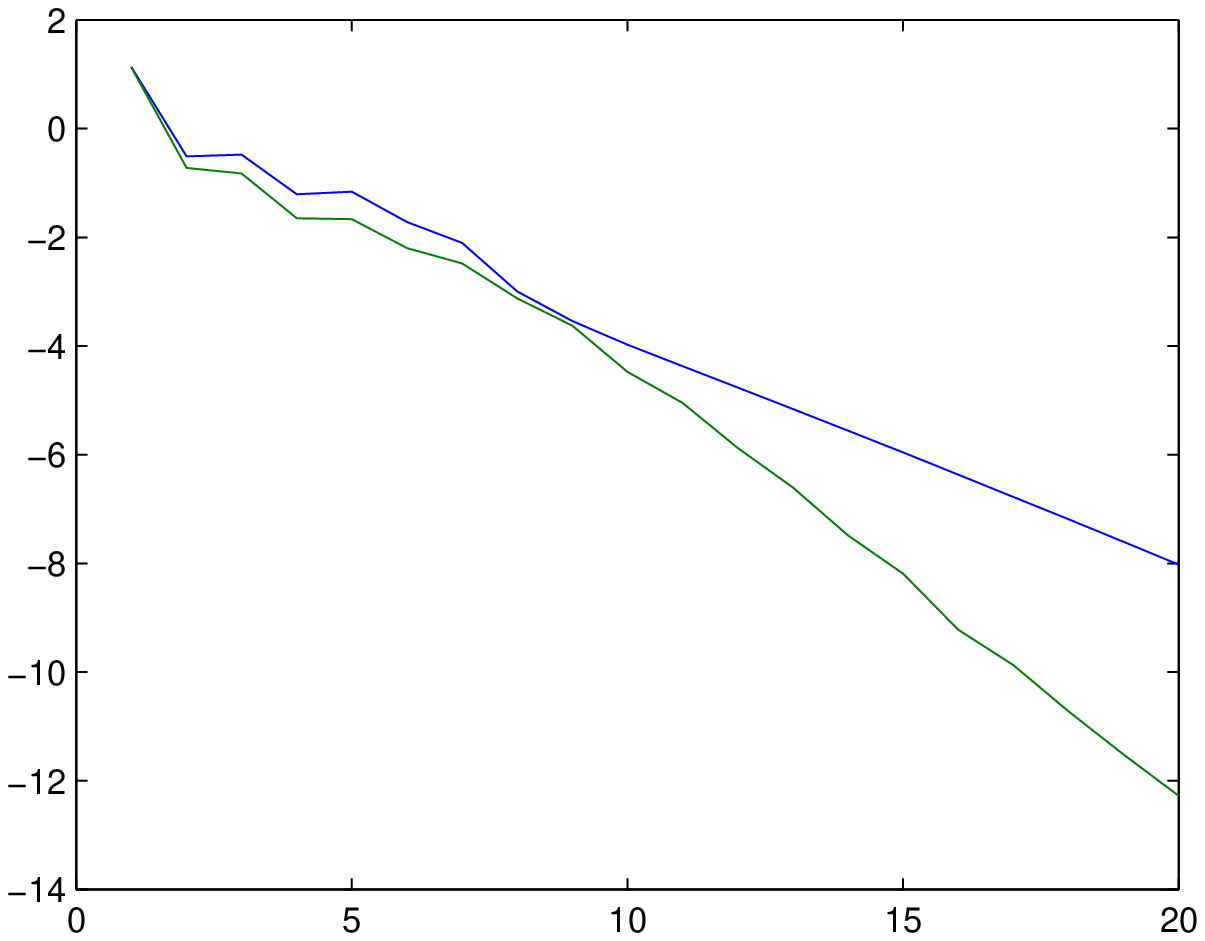}
\caption{\label{fig:eps-3} Log of the error on the computed solution as a function of the number of Schwarz iterations for Rossby number $\eps=10^{-3}$ (left)  and $\eps=10^{-2}$ (right) and for a random initial guess using the Taylor approximation parameter $\alpha_{Tay}$ (up) and the optimized one $\alpha_{opt}$ (down)} 
\end{figure}

\begin{figure}[hbtp]
\centering
\includegraphics[height=.35\hsize,angle=0]{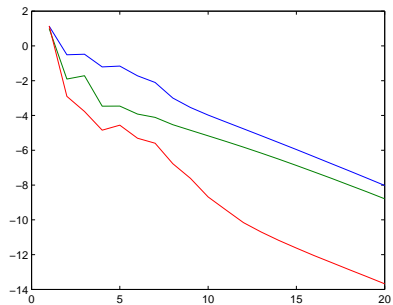}
\caption{\label{fig:sin}Log of the error on the computed solution 
as a function of the number of Schwarz iterations 
for Rossby number $\eps=10^{-2}$ and with optimal parameter for a low frequency signal (down)
and for a high frequency signal (middle) and for a random signal (up)}
\end{figure}

In this section we  are interested 
in the optimization the transmission conditions (\ref{Bpm})
with respect to the free parameters $\alpha$ and $\beta$.
To optimize the conditions means that we choose parameters $\alpha$ and $\beta$
such that the Schwarz waveform relaxation algorithm (\ref{Algo}) reaches 
a given error for as small as possible number of iterations.
The analytical solution of this problem is quite complex in the considered framework
and we only present here a numerical strategy to reach the optimum.
In the simpler case of a 1D advection diffusion equation
a complete solution of the related optimization problem is given in \cite{halpern07}.\\ 
\newline
We consider a test case for which all the initial data (velocities and perturbation of the water height)
are taken equal to zero. We initialize the algorithm (\ref{swra}) with random boundary conditions on the interface
and we study the convergence of the solution towards the analytical ones.
This test is quite classical to study the convergence of a domain decomposition algorithm.
It is interesting since the initial quantities do contain all frequencies.
In all the computations the physical parameters $Re$ and $Fr$ are taken equal to one
but the Rossby number $\eps$ remains free. For a given value of $\eps$
we apply the transmission conditions (\ref{Bpm})
for several values of the parameters $\alpha$ and $\beta$ and we compare the $L^2$ error 
between the computed and the analytical solutions after a given number of iterations.
It allows us to find an optimal pair $(\alpha_{opt},\beta_{opt})$ that minimizes this error.
This first study exhibit that the influence of the parameter $\beta$ is quite small.
In the following this parameter will be kept equal to its theoretical value (\ref{defab}).
In a second step we study the dependency of the optimal parameter $\alpha_{opt}$ 
with respect to the Rossby number $\eps$.
The results are presented in Fig. \ref{fig:alphaepsilon} for different values of $\eps$. We found that this optimized parameter does depend on $\eps$ in a nontrivial way. 
\newline
We now present the evolution of the error on the computed solution as a function of the number of iterations of the Schwarz waveform relaxation algorithm (\ref{swra}) in both cases $\alpha=\alpha_{opt}$ and $\alpha=\alpha_{Tay}$. 
in Fig. \ref{fig:eps-3} we present the results 
for two different values of the Rossby number $\eps$ :
$\eps=10^{-3}$ and $\eps=10^{-2}$.
The curves (Log of the error) all look like straight lines, 
at least after a sufficiently large number of iterations.
The method appears to be more efficient when the Rossby number is smaller 
since the error decreases much faster in the case $\eps=10^{-3}$ 
- Fig. \ref{fig:eps-3} on the left.
This result is consistent with the previous theoretical study 
that is based on an asymptotic analysis in $\eps$.
We also observe that for a given value of $\eps$ 
the curves look similar for both optimized and Taylor approximation parameters
even if the error decreases faster for the optimal value $\alpha_{opt}$. 
Moreover let us observe that to reach an error of $10^{-4}$
(that is enough for the applicability of the Schwarz waveform relaxation algorithm)
both algorithms (with optimized or Taylor approximation parameter) 
need a very close number of iterations.\\
\newline
We compute the same test with Rossby number $\eps=10^{-2}$ 
but with a sinusoidal initial guess (instead of the random ones) 
for the transmission conditions.
We consider two different sinusoids with one or ten periods 
in the space-time considered interval
and we use Taylor approximation parameters $\alpha_{Tay}$ and $\beta_{Tay}$.
In Fig. \ref{fig:sin} the results appears to be much better 
for the low frequency sinusoid as for high frequency one.
The results for the high frequency sinusoid look similar 
to the results that were obtained with the random initial guess.
It follows that the method is particularly well adapted to low frequency signals :
the relative error is smaller than $10^{-4}$ after only two iterations.
\subsection{Numerical application}
\begin{figure}[hbtp]
\centering
\includegraphics[width=.75\textwidth,angle=-0]{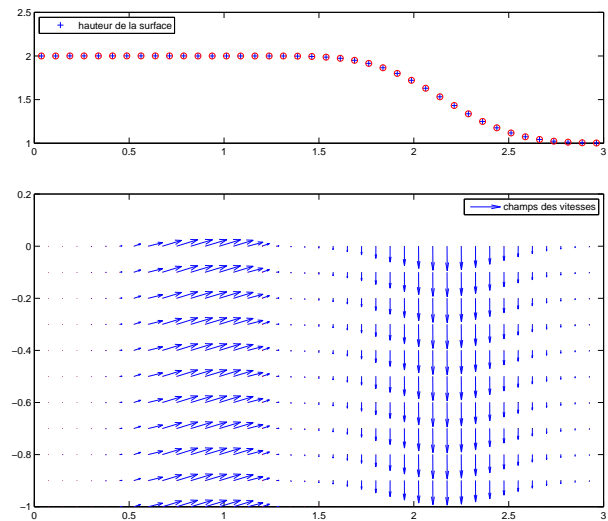}
\caption{\label{fig:final}Water height and velocity field at final time}
\end{figure}

\begin{figure}[hbtp]
\centering
\includegraphics[width=.5\textwidth,angle=-0]{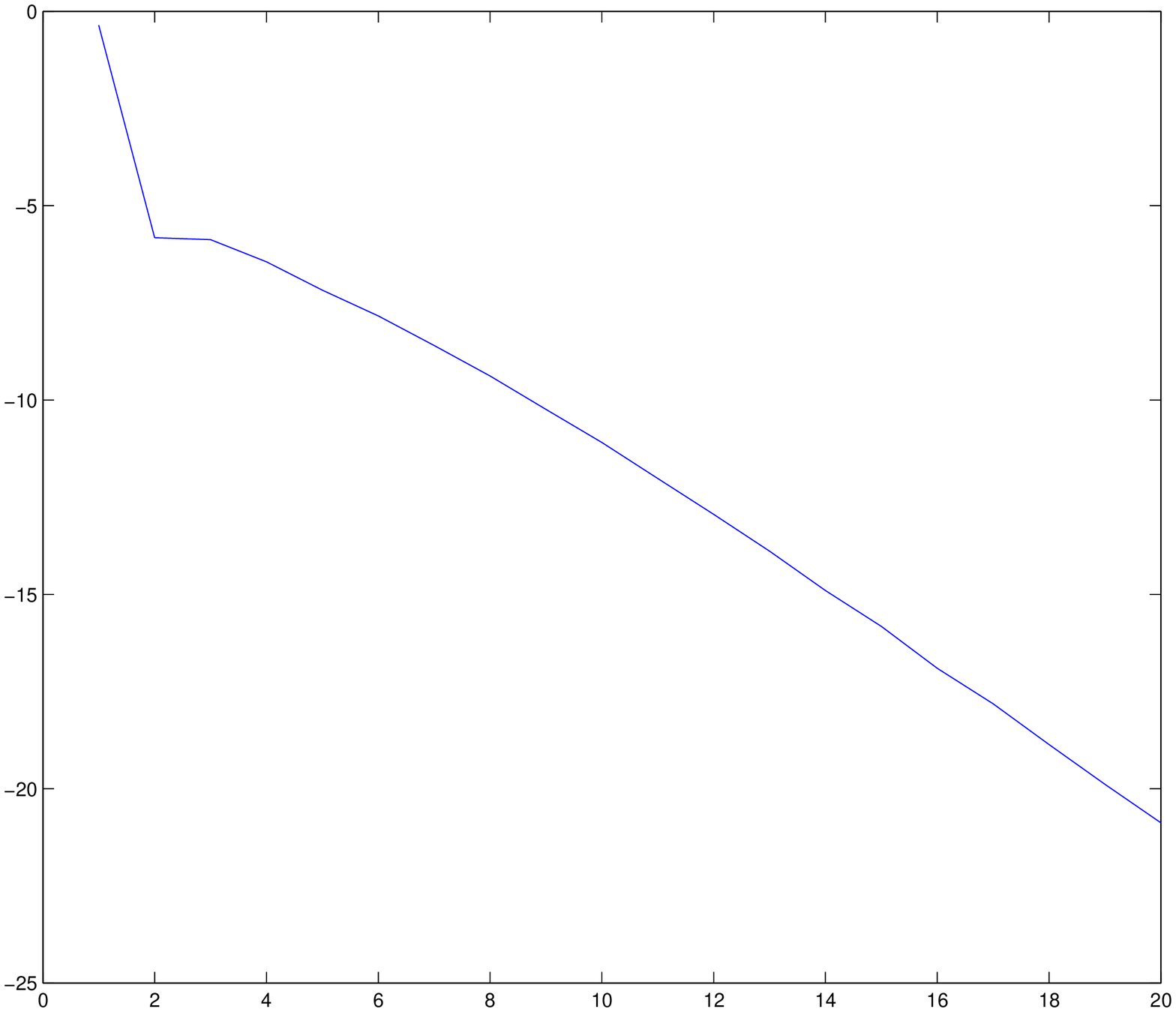}
\caption{Log of the relative error on the solution computed by using the Schwarz waveform relaxation algorithm versus number of iterations}
\label{fig:herr}
\end{figure}

In this section we consider the case of a flow 
with a constant positive background velocity 
$u_0=1. m/s$ and an initial local decreasing step on the water height. 
We choose the Rossby number $\eps$ equal to $10^{-3}$.
We choose $nx=40$, $nz=10$ and $nt=40$ in order to ensure the CFL condition.
We present the initial solution 
and the solution computed at final time $T=1.3 s$ 
after 20 iterations by the proposed Schwarz waveform relaxation algorithm
in Fig. \ref{fig:final}. The 2d horizontal velocity vector field $(u,v)$ is presented 
in the 2d vertical domain (in the $(x,z)$ plane) which is occupied by the flow. 
A horizontal vector denotes a velocity which is collinear to the $x$-direction 
and a vertical one denotes a velocity which is collinear to the $y$-direction.
Since we consider the linearized version of the equations
the step just moves without deformation from the left to the right of the domain.
Since the Coriolis effect is dominant we observe the formation 
of a transverse jet which moves with the step.
Another consequence of the Coriolis effect is the formation of a stationary eddy
at the initial location of the step.
\begin{figure}[hbtp]
\centering
\includegraphics[width=.75\textwidth,angle=-0]{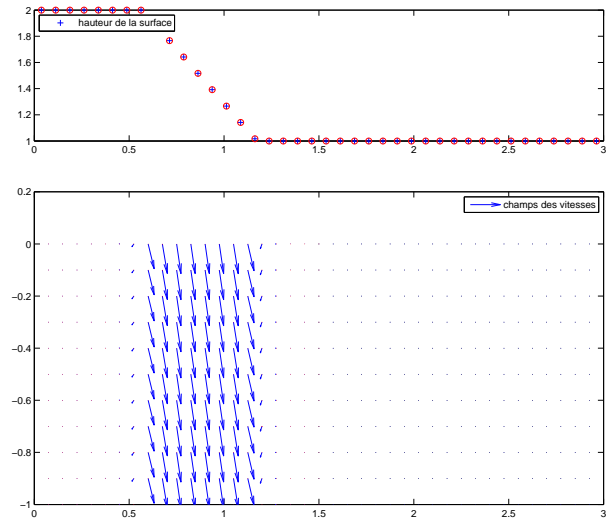}
\caption{\label{fig:initial}Water height and velocity field at initial time}
\end{figure}
We now compare the solution that is computed on the whole domain 
with the solution that is obtained by considering 
the presented domain decomposition strategy.
In Fig. \ref{fig:herr} we present the evolution of the relative error 
between the two solutions versus the number of considered iterations.
It exhibits the fast convergence of the algorithm for such a case.
After two iterations the relative error is around $10^{-6}$ 
and it reaches the factor $10^{-10}$ after eight iterations.
\section{Conclusion}
We presented in this article a new domain decomposition method 
for the viscous primitive equations. 
It involves a Schwarz waveform relaxation type algorithm 
with approximated transmission conditions
for which we proved well-posedness.
We presented a numerical optimization of the transmission conditions
and we study the speed of convergence of the algorithm for several test cases.
Academic numerical applications were presented.
In forthcoming papers we plan to 
prove the convergence of the algorithm and 
we want to present oceanographic configurations
and to increase the efficiency of the algorithm
by deriving more complex transmission conditions
based on another asymptotic regime
that corresponds to quasi-geostrophic flows.
\label{secconc}
\section*{Acknowledgements}
The authors thank L. Halpern and V. Martin for fruitful discussions and helpful comments. This work was partially supported by ANR program COMMA (http://www-lmc.imag.fr/COMMA/). 


\begin{thebibliography}{99}


\bibitem{arakawa}
Arakawa A. \& Lamb V.
Computational design of the basic dynamical processes of the UCLA general circulation model,
{\it Methods in Computational Physics}, Vol. 17 (1977), pp 174--267.

\bibitem{roms}
Arango H.G. \& Shchepetkin A.F.,
ROMS : A Regional Ocean Modeling System. 
[www.myroms.org/index.php]

\bibitem{blayo}
Cailleau S., Fedorenko V., Barnier B., Blayo E. \& Debreu L.  Comparison of different numerical methods used to handle the open boundary of a regional ocean circulation model of the Bay of Biscay, \textit{submitted}  (2007) 


\bibitem{titi}
Cao C. \& Titi E.S., 
Global well-posedness of the three-dimensional viscous primitive equations 
of large scale ocean and atmosphere dynamics, 
{\it Annals of Mathematics}, Vol. 166 (2007), No. 1, pp 245--267.

\bibitem{C-R94}
Cushman-Roisin, B. {\it Introduction to Geophysical Fluid Dynamics}, Prentice Hall (1994), pp 320.

\bibitem{gander01}
Daoud D.S. \& Gander M.J.,
Overlapping Schwarz waveform relaxation for convection reaction diffusion problems,
{\it Proceedings of the 13th International Conference on Domain Decomposition Methods},
2001, pp 253--260. 
[www.ddm.org/conferences.html]

\bibitem{engquist}
Engquist B. \& Majda A.,
{Absorbing boundary conditions for the numerical simulation of waves},
{\it Math. Comp.}, Vol. 31 (1977), No. 139, pp 629--651.

\bibitem{gander97}
Gander M.J.,
Overlapping Schwarz for parabolic problems,
{\it Proceedings of the 9th International Conference on Domain Decomposition Methods},
(1997), pp 97--104. 
[www.ddm.org/conferences.html]

\bibitem{gander98}
Gander M.J.,
{A waveform relaxation algorithm with overlapping splitting for reaction diffusion equations},
{\it Numerical Linear Algebra with Applications}, Vol. 6 (1998), pp 125--145.

\bibitem{gander06}
Gander M.J.,
{Optimized Schwarz methods},
{\it SIAM Journal of Numerical Analysis}, Vol. (2006), No., pp 699--731.

\bibitem{halpern07}
Gander M.J. \& Halpern L.,
Optimized Schwarz Waveform Relaxation for Advection Reaction Diffusion Problems,  
{\it SIAM Journal on Numerical Analysis}, Vol. 45 (2007), No. 2, pp 666--697.

\bibitem{halpern03}
Gander M.J., Halpern L. \& Nataf F.,
{Optimal Schwarz waveform relaxation for the one dimensional wave equation},
{\it SIAM Journal of Numerical Analysis}, Vol. 41 (2003), No. 5, pp 1643--1681.

\bibitem{perthame} 
Gerbeau J.F. \& Perthame B.,
{Derivation of viscous Saint Venant system for
laminar shallow water; numerical simulation}, 
{\it Discrete and Continuous Dynamical Systems - Series B}, Vol. 1 (2001), No. 1, pp 89--102.

\bibitem{giladi}
Giladi E. \& Keller H.B.,
{Space time domain decomposition for parabolic problems},
{\it Numerische Mathematik}, Vol. 93 (2002), No. 2, pp 279--313. 

\bibitem{halpern86}
Halpern L.,
{Artificial boundary conditions for the advection diffusion equations},
{\it Math. Comp.}, Vol. 174 (1986), pp 425--438.

\bibitem{halpern91}
Halpern L., 
Artificial boundary conditions for incompletely parabolic perturbations of hyperbolic systems, 
{\it SIAM Journal on Math. Anal.}, Vol. 22 (1991), No. 5, pp 1256--1283.

\bibitem{jeltsch}
Jeltsch R. \& Pohl B.,
Waveform relaxation with overlapping splittings,
{\it SIAM J. Sci. Comp.}, Vol. 16 (1995), No. 1, pp 40--49.

\bibitem{lelarasmee}
Lelarasmee E., Ruehli A.E. \& Sangiovanni Vincetelli A.L.,
{The waveform relaxation method for time-domain analysis 
of large scale integrated circuits},
{\it IEEE Trans. on CAD of IC and Systems}, Vol. 1 (1982), pp 131--145.

\bibitem{lions1}
Lions P.L.,
{On the Schwarz alternating method I},
{\it Chan T.F., Glowinski R., Periaux J \& Widlund O. editors,
Proceedings of the 1st International Conference on Domain Decomposition Methods},
SIAM, (1988).

\bibitem{lions3}
Lions P.L.,
{On the Schwarz alternating method II: 
a variant for nonoverlapping subdomains},
{\it Chan T.F., Glowinski R., Periaux J \& Widlund O. editors,
Proceedings of the 3rd International Conference on Domain Decomposition Methods},
SIAM, (1990).

\bibitem{lionspere}
Lions J.L., Temam R. \& Wang S.
New formulations of the primitive equations of the atmosphere and applications,
{\it Nonlinearity}, Vol. 5 (1992), pp 237--288.

\bibitem{rousseau}
Lucas C. \& Rousseau A.,
New developments and cosine effect in the viscous shallow water and quasi geostrophic equations,
{\it submitted}.

\bibitem{nemo}
Madec G., Delecluse P., Imbard M. \& Lévy C., 
{\it 1998: OPA 8.1 Ocean General Circulation Model reference manual}, 
{Note du Pole de mod\'elisation, Institut Pierre-Simon Laplace (IPSL), France}, No. 11, 91pp, 1998.
[www.locean-ipsl.upmc.fr/NEMO].

\bibitem{martin05}
Martin V.,
{An optimized Schwarz waveform relaxation method for unsteady 
convection diffusion equation},
{\it Applied Numerical Mathematics}, Vol. 52 (2005), No. 4, pp 401--428.

\bibitem{martin08}
 Martin V., {A Schwarz Waveform Relaxation Method for the Viscous Shallow Water Equations}, 
{\it Domain Decomposition Methods in Science and Engineering}, Vol. 40 (2004), pp 653--660. 

\bibitem{mom}
Pacanowski R.C. \& Griffies S.M.,
{\it MOM 3.0 Manual}, (2000).
[www.gfdl.noaa.gov/~smg/MOM/web/guide\_parent].

\bibitem{quarteroni}
Quarteroni A. \& Valli A.,
{\it Domain Decomposition Methods for PDEs},
Oxford Science Publications, London, (1999).

\bibitem{stvenant}
de Saint-Venant A.J.C., 
Th{\'e}orie du mouvement non-permanent des eaux, avec application 
aux crues des rivi{\`e}res et \`a l'introduction des mar{\'e}es dans leur lit (in french), 
{\it C. R. Acad. Sc., Paris},
Vol. 73 (1871), pp 147--154.

\bibitem{schwarz}
Schwarz H.A.,
{\"Uber einen Grenz\"ubergang durch alternierendes Verfahren},
{\it Vierteljahrschrift der Naturforschenden Gesellschaft in Z\"urich},
Vol. 15 (1870), pp 272--286.


\bibitem{SW97}
Showalter, R. E., {\it Monotone operators in {B}anach space and nonlinear partial differential equations}, {Mathematical Surveys and Monographs} Vol 49. (1997) pp xiv+278

\bibitem{temam08}
Temam R. \& Tribbia J., 
{\it Computational methods for the oceans and the atmosphere},
Ciarlet P.G. General Editor,
Special volume of the Handbook of numerical analysis,
Elsevier, Amsterdam, (2008).

\bibitem{temam04}
Temam R. \& Ziane M.,
{\it Some mathematical problems in geophysical fluid dynamics},
Friedlander S. \& Serre D. editors,
Handbook of Mathematical Fluid Dynamics, Vol. 3, Elsevier, (2004).

\bibitem{toselli}
Toselli A. \& Widlund O.,
{\it Domain decomposition methods - Algorithms and theory},
Series in Computational Mathematics, Vol. 34, Springer, (2004).

\end{thebibliography}
\end{document}